\documentstyle[11pt]{article}

\input amssym 
\catcode`\@=11\relax
\newwrite\@unused
\def\typeout#1{{\let\protect\string\immediate\write\@unused{#1}}}

\def\psglobal#1{
\immediate\special{ps: plotfile #1 }}
\def\psfiginit{\typeout{psfiginit}
\immediate\psglobal{figtex.pro}%
\special{ps:: /TeXMagnification {\the\mag} def}
}

\def\@nnil{\@nil}
\def\@empty{}
\def\@psdonoop#1\@@#2#3{}
\def\@psdo#1:=#2\do#3{\edef\@psdotmp{#2}\ifx\@psdotmp\@empty \else
    \expandafter\@psdoloop#2,\@nil,\@nil\@@#1{#3}\fi}
\def\@psdoloop#1,#2,#3\@@#4#5{\def#4{#1}\ifx #4\@nnil \else
       #5\def#4{#2}\ifx #4\@nnil \else#5\@ipsdoloop #3\@@#4{#5}\fi\fi}
\def\@ipsdoloop#1,#2\@@#3#4{\def#3{#1}\ifx #3\@nnil
       \let\@nextwhile=\@psdonoop \else
      #4\relax\let\@nextwhile=\@ipsdoloop\fi\@nextwhile#2\@@#3{#4}}
\def\@tpsdo#1:=#2\do#3{\xdef\@psdotmp{#2}\ifx\@psdotmp\@empty \else
    \@tpsdoloop#2\@nil\@nil\@@#1{#3}\fi}
\def\@tpsdoloop#1#2\@@#3#4{\def#3{#1}\ifx #3\@nnil
       \let\@nextwhile=\@psdonoop \else
      #4\relax\let\@nextwhile=\@tpsdoloop\fi\@nextwhile#2\@@#3{#4}}
\def\psdraft{
	\def\@psdraft{0}
	\def\@psdraftspecial{100}
}
\def\psdraftspecial{
	\def\@psdraft{0}
	\def\@psdraftspecial{0}
}
\def\psfull{
	\def\@psdraft{100}
}
\psfull

\newif\if@prologfile
\newif\if@postlogfile
\newif\if@bbllx
\newif\if@bblly
\newif\if@bburx
\newif\if@bbury
\newif\if@height
\newif\if@width
\newif\if@rheight
\newif\if@rwidth
\newif\if@clip
\newif\if@right
\newif\if@left
\newif\if@toplines
\newif\if@box
\newif\if@caption
\newif\if@surround
\newif\if@captionwidth
\newif\if@captionwrite
\newif\if@captionopen
\def\@p@@sclip#1{\@cliptrue}
\def\@p@@sfile#1{
		\def\@p@sfile{#1}
}
\def\@p@@sfigure#1{
		\def\@p@sfile{#1}
}
\def\@p@sfake{\hbox to 0pt{\hss Whatever\hss}}
\def\@p@@sbbllx#1{
		\@bbllxtrue
		\@d@mscratch=#1
		\edef\@p@sbbllx{\number\@d@mscratch}
}
\def\@p@@sbblly#1{
		\@bbllytrue
		\@d@mscratch=#1
		\edef\@p@sbblly{\number\@d@mscratch}
}
\def\@p@@sbburx#1{
		\@bburxtrue
		\@d@mscratch=#1
		\edef\@p@sbburx{\number\@d@mscratch}
}
\def\@p@@sbbury#1{
		\@bburytrue
		\@d@mscratch=#1
		\edef\@p@sbbury{\number\@d@mscratch}
}
\def\@p@@sheight#1{
		\@heighttrue
		\@d@mscratch=#1
   		\edef\@p@sheight{\number\@d@mscratch}
}
\def\@p@@swidth#1{
		\@widthtrue
		\@d@mscratch=#1
		\edef\@p@swidth{\number\@d@mscratch}
}
\def\@p@@srheight#1{
		\@rheighttrue
		\@d@mscratch=#1
		\edef\@p@srheight{\number\@d@mscratch}
}
\def\@p@@srwidth#1{
		\@rwidthtrue
		\@d@mscratch=#1
		\edef\@p@srwidth{\number\@d@mscratch}
}
\def\@p@@sright#1{\@righttrue \@surroundtrue}
\def\@p@@sleft#1{\@lefttrue \@surroundtrue}
\def\@p@@sextraheight#1{\@d@mextraheight=#1}
\def\@p@@sbox#1{\@boxtrue}
\def\@p@@scaption#1{\@captiontrue}
\def\@p@@stoplines#1{
		\@toplinestrue
		\@c@ttoplines=#1
}
\def\@p@@scaptionwidth#1{
		\@captionwidthtrue
	  	\@d@mcaptionwidth=#1
}
\def\@p@@scaptionwrite#1{
		\global\@captionwritetrue
		\global\@w@rname=\expandafter{\jobname_captions.tex}
		\typeout{Captions are written to \the\@w@rname}
}

\def\@p@@sprolog#1{\@prologfiletrue\def\@prologfileval{#1}}
\def\@p@@spostlog#1{\@postlogfiletrue\def\@postlogfileval{#1}}
\def\@cs@name#1{\csname #1\endcsname}
\def\@setparms#1=#2,{\@cs@name{@p@@s#1}{#2}}
\def\ps@init@parms{
		\@bbllxfalse \@bbllyfalse
		\@bburxfalse \@bburyfalse
		\@heightfalse \@widthfalse
		\@rheightfalse \@rwidthfalse
		\def\@p@sbbllx{}\def\@p@sbblly{}
		\def\@p@sbburx{}\def\@p@sbbury{}
		\def\@p@sheight{}\def\@p@swidth{}
		\def\@p@srheight{}\def\@p@srwidth{}
		\def\@p@sfile{}
		\def\@p@scost{10}
		\def\@sc{}
		\@prologfilefalse
		\@postlogfilefalse
		\@clipfalse
		\@rightfalse \@leftfalse
		\@boxfalse \@captionfalse
		\@toplinesfalse \@surroundfalse
		\@d@mextraheight=0pt
 		\@c@ttoplines=0
		\@pshape={} \def\@p@srheight@total{}
		\@captionwidthfalse \@d@mcaptionwidth=0pt
}
\def\parse@ps@parms#1{
	 	\@psdo\@psfiga:=#1\do
		   {\expandafter\@setparms\@psfiga,}}
\newif\ifno@bb
\newif\ifnot@eof
\newread\ps@stream
\newtoks\@linetok
\def\bb@missing{
	\typeout{psfig: searching \@p@sfile \space  for bounding box}
	\openin\ps@stream=\@p@sfile
	\no@bbtrue
	\not@eoftrue
	\catcode`\%=12
	\loop
		\read\ps@stream to \line@in
		\global\@linetok=\expandafter{\line@in}
		\ifeof\ps@stream \not@eoffalse \fi
		\@bbtest{\@linetok}
		\if@bbmatch\not@eoffalse\expandafter\bb@cull\the\@linetok\fi
	\ifnot@eof \repeat
	\catcode`\%=14
}	
\catcode`\%=12
\newif\if@bbmatch
\def\@bbtest#1{\expandafter\@a@\the#1
\long\def\@a@#1
     \ifx\@bbtest#2\@bbmatchfalse\else\@bbmatchtrue\fi}
\long\def\bb@cull#1 #2 #3 #4 #5 {
	\@d@mscratch=#2 bp\edef\@p@sbbllx{\number\@d@mscratch}
	\@d@mscratch=#3 bp\edef\@p@sbblly{\number\@d@mscratch}
	\@d@mscratch=#4 bp\edef\@p@sbburx{\number\@d@mscratch}
	\@d@mscratch=#5 bp\edef\@p@sbbury{\number\@d@mscratch}
	\no@bbfalse
}
\catcode`\%=14
\def\compute@bb{
		\no@bbfalse
		\if@bbllx \else \no@bbtrue \fi
		\if@bblly \else \no@bbtrue \fi
		\if@bburx \else \no@bbtrue \fi
		\if@bbury \else \no@bbtrue \fi
		\ifno@bb \bb@missing \fi
		\ifno@bb \typeout{FATAL ERROR: no bb supplied or found}
			\no-bb-error
		\fi
		\count203=\@p@sbburx
		\count204=\@p@sbbury
		\advance\count203 by -\@p@sbbllx
		\advance\count204 by -\@p@sbblly
		\edef\@bbw{\number\count203}
		\edef\@bbh{\number\count204}
}
\def\in@hundreds#1#2#3{\count240=#2 \count241=#3
		     \count100=\count240	
		     \divide\count100 by \count241
		     \count101=\count100
		     \multiply\count101 by \count241
		     \advance\count240 by -\count101
		     \multiply\count240 by 10
		     \count101=\count240	
		     \divide\count101 by \count241
		     \count102=\count101
		     \multiply\count102 by \count241
		     \advance\count240 by -\count102
		     \multiply\count240 by 10
		     \count102=\count240	
		     \divide\count102 by \count241
		     \count200=#1\count205=0
		     \count201=\count200
			\multiply\count201 by \count100
		     	\advance\count205 by \count201
		     \count201=\count200
			\divide\count201 by 10
		     	\multiply\count201 by \count101
			\advance\count205 by \count201
		     \count201=\count200
			\divide\count201 by 100
			\multiply\count201 by \count102
			\advance\count205 by \count201
		     \edef\@result{\number\count205}
}
\def\compute@wfromh{
		\in@hundreds{\@p@sheight}{\@bbw}{\@bbh}
		\edef\@p@swidth{\@result}
}
\def\compute@hfromw{
		\in@hundreds{\@p@swidth}{\@bbh}{\@bbw}
		\edef\@p@sheight{\@result}
}
\def\compute@handw{
		\if@height
			\if@width
			\else
				\compute@wfromh
			\fi
		\else
			\if@width
				\compute@hfromw
			\else
				\edef\@p@sheight{\@bbh}
				\edef\@p@swidth{\@bbw}
			\fi
		\fi
}
\def\compute@resv{
		\if@rheight \else \edef\@p@srheight{\@p@sheight} \fi
		\if@rwidth \else \edef\@p@srwidth{\@p@swidth} \fi
		\edef\@p@srheight@total{\@p@srheight}
}
\newtoks\@pshape
\def\@c@ttoplines{\count120}
\def\@c@theightcount{\count121}
\def\@c@tshapecount{\count122}
\newdimen\@d@mwidthshape
\newdimen\@d@mextraheight
\newdimen\@d@mscratch
\def\compute@parshape{
	\if@right
		\if@left
	   		\typeout{error: Can't have both left and right set}
			\@leftfalse
		\fi
	\fi
	\@d@mscratch=\@p@swidth truesp
	\divide \@d@mscratch by 19
	\multiply \@d@mscratch by 20
	\edef\@p@swidthdimen{\the\@d@mscratch}
	\@c@tshapecount=\@c@ttoplines
 	\@d@mscratch=\baselineskip
	\multiply \@d@mscratch by \@c@ttoplines
	\advance \@d@mscratch by .4\baselineskip
    	\edef\@p@stopdistance{\the\@d@mscratch }
	\@d@mscratch=\@p@sheight truesp
	\divide \@d@mscratch by 2
	\edef\@p@shalfboxheight{\the\@d@mscratch}
	\if@toplines
		\loop \@pshape=\expandafter{\the\@pshape 0pt \hsize}
		\advance\@c@ttoplines by -1
		\ifnum\@c@ttoplines>0 \repeat
	\fi
   	\@c@theightcount=\@p@srheight@total
	\advance \@c@theightcount by \@d@mextraheight
	\divide  \@c@theightcount by \baselineskip
	\advance \@c@theightcount by 1
    	\advance \@c@tshapecount by \@c@theightcount
	\advance \@c@theightcount by -1
	\@d@mwidthshape=\hsize
     	\advance \@d@mwidthshape by -\@p@swidthdimen
	\if@left
		\def\@moveshape{0pt}
		\ifnum\@c@theightcount>0
		      	\loop
			\@pshape=%
			\expandafter{\the\@pshape %
					\@p@swidthdimen \@d@mwidthshape}
			\advance \@c@theightcount by -1
			\ifnum\@c@theightcount>0 \repeat
		\else
			\advance \@c@tshapecount by 1
		\fi
	\fi
	\if@right
		\@d@mscratch=\hsize
		\advance \@d@mscratch by -\@p@swidth truesp
		\edef\@moveshape{\@d@mscratch}
		\ifnum\@c@theightcount>0
			\loop
			\@pshape=\expandafter{\the\@pshape 0pt \@d@mwidthshape}
			\advance \@c@theightcount by -1
			\ifnum\@c@theightcount>0 \repeat
		\else
			\advance \@c@tshapecount by 1
		\fi
	\fi
	\ifnum \@p@srheight=0
		\@pshape={}
		\@c@tshapecount = 0
	\else
	 	\@pshape=\expandafter{\the\@pshape 0pt \hsize}
	\fi
}
\def\@p@ssetsurroundboxes{
	\global\parshape=\count122 \the\@pshape		
 	\moveright\@moveshape
	\vbox to 0pt\bgroup\hskip0pt\vskip\@p@stopdistance
}
\newtoks\@captiontok
\newbox\@b@xcaption
\newdimen\@d@mcaptionwidth
\newdimen\@d@mcaptionheight
\newwrite\@w@rcaption
\newtoks\@w@rname
\def\setcaption#1{\@captiontok={#1}}
\def\@set@caption{
	\setbox\@b@xcaption\vbox{\hsize\@d@mcaptionwidth
	\tolerance=9000 \baselineskip=11.4pt
	\noindent\relax\the\@captiontok}
	\if@captionwrite
		\if@captionopen
		\else
			\global\@captionopentrue
			\immediate\openout\@w@rcaption=\the\@w@rname
		\fi
		\immediate\write\@w@rcaption{\the\@captiontok}
		\immediate\write\@w@rcaption{}
	\fi
}
\def\compute@caption{
	\if@captionwidth
	\else
		\@d@mcaptionwidth = \@p@swidth truesp
		\divide \@d@mcaptionwidth by 20
		\multiply \@d@mcaptionwidth by 17
	\fi
	\@set@caption
	\@d@mcaptionheight=\ht\@b@xcaption
	\if@rheight
	\else
		\count100=\@p@srheight
	   	\advance \count100 by \@d@mcaptionheight
	   	\advance \count100 by \bigskipamount
	   	\advance \count100 by \medskipamount
	   	\edef\@p@srheight@total{\number\count100}
	\fi
}
\newif\if@alreadyjtem \@alreadyjtemfalse
\def\newpar{\hfil\vadjust{\vskip\parskip}%
	{\count100=\parskip
	\count101=\baselineskip
	\divide\count101 by 10  \multiply\count101 by 3
	\advance \count100 by \count101
	\divide\count100 by \baselineskip
	\advance\count100 by \prevgraf
	\global\prevgraf=\count100}%
	\break\if@alreadyjtem\else\indent\fi%
}
\let\sav@par=\par
\def\jtem#1{%
    	\if@alreadyjtem\else\bgroup\fi
	\def\par{\sav@par\egroup\sav@par}
	\if@alreadyjtem\else\leftskip\parindent\fi
	\@alreadyjtemtrue
	\noindent\hskip0pt
	\llap{#1\ }\ignorespaces
}
\def\compute@sizes{%
	\compute@bb
	\compute@handw
  	\compute@resv
	\if@caption
		\compute@caption
	\fi
	\if@surround
		\compute@parshape
	\fi
}
\def\@p@sdospecials{
	\ifnum\@p@scost<\@psdraft
	       	\typeout{psfig: including \@p@sfile \space }
	\fi
	\special{ps::[begin] 	\@p@swidth \space \@p@sheight \space
			\@p@sbbllx \space \@p@sbblly \space
			\@p@sbburx \space \@p@sbbury \space
			startTexFig \space }
	\ifnum\@p@scost<\@psdraft
		\if@clip
			\typeout{(clip)}
			\special{ps:: \@p@sbbllx \space \@p@sbblly \space
				\@p@sbburx \space \@p@sbbury \space
			    	doclip \space }
		\fi
	\fi
	\if@box
		\typeout{(box)}
  		\special{ps:: \@p@sbbllx \space \@p@sbblly \space
			\@p@sbburx \space \@p@sbbury \space
		    	dobox \space }
	\fi
	\ifnum\@p@scost<\@psdraft
		\if@prologfile
	    		\special{ps: plotfile \@prologfileval \space }
		\fi
		\special{ps: plotfile \@p@sfile \space }
    		\if@postlogfile
			\special{ps: plotfile \@postlogfileval \space }
		\fi
	\fi
	\special{ps::[end] endTexFig \space }
}
\newif\if@putinvbox

\def\psfig#1{{%
	\ifhmode%
		\vbox\bgroup
		\@putinvboxtrue
	\else
		\@putinvboxfalse
	\fi
       	\ps@init@parms
	\parse@ps@parms{#1}
       	\compute@sizes
	\if@surround
		\psfig@for@surround
	\else
		\psfig@for@regular
	\fi
	\if@putinvbox
       		\egroup
	\fi
}}
\def\psfig@for@surround{%
	\@p@ssetsurroundboxes
	\ifnum\@p@scost<\@psdraft
		\@p@sdospecials
		\vbox to \@p@srheight true sp{\vss}
       	\else
		\if@box
			\@p@sdospecials
		\fi
		\vbox to \@p@srheight true sp{
			\vskip\@p@shalfboxheight
			\hbox to \@p@srwidth true sp{
				\hss
				\ifnum\@p@scost<\@psdraftspecial
					\@p@sfile
				\else
					\@p@sfake
				\fi
      				\hss
			}
		\vss
		}
	\fi
	\if@caption
		\medskip
		\hbox to \@p@srwidth true sp{
			\hss
			\box\@b@xcaption
			\hss
		}
 		\medskip
	\fi
	\vss\egroup
	\vskip-\parskip
}

\def\psfig@for@regular{%
	\if@putinvbox
	\else
		\vskip\parskip
	\fi
	\ifnum\@p@scost<\@psdraft
		\@p@sdospecials
		\vbox to \@p@srheight true sp{%
			\hbox to \@p@srwidth true sp{
			\hfil
			}
		\vfil
		}
       	\else
		\if@box
			\@p@sdospecials
		\fi
	    	\vbox to \@p@srheight true sp{
			\vss
			\hbox to \@p@srwidth true sp{
				\hss
				\ifnum\@p@scost<\@psdraftspecial
					\@p@sfile
				\else
					\@p@sfake
				\fi
				\hss
			}
		    	\vss
		}
	\fi
	\if@caption
		\medskip
		\hbox to \@p@srwidth true sp{
			\hss
			\box\@b@xcaption
			\hss
		}
		\bigskip
	\fi
	\if@putinvbox
	\else
		\vskip-\parskip
	\fi
}
\catcode`\@=12\relax
\psfiginit

\font\tinybbfont=msbm6
\font\scriptsizebbfont=msbm7 scaled \magstep 1
 1
\font\footnotesizebbfont=msbm9 scaled \magstep 0
 2
\font\bbfont=msbm9 scaled \magstep1  
 1
 2
 3
 4

\def\tinyBbb#1{\hbox{\tinybbfont #1}}
\def\scriptsizeBbb#1{\hbox{\scriptsizebbfont #1}}

\def\footnotesizeBbb#1{\hbox{\footnotesizebbfont #1}}

\def\Bbb#1{\hbox{\bbfont #1}}

\newcommand{\CP}{{\Bbb C}{\rm P}}

\newcommand{\GL}{\mbox{\it GL}\,}
\newcommand{\Hom}{\mbox{\it Hom}\,}

\newcommand{\Pic}{\mbox{\rm Pic}\,}
\newcommand{\Proj}{\mbox{\rm Proj}\,}

\newcommand{\SL}{\mbox{\it SL}}

\newcommand{\Stab}{\mbox{\it Stab}\,}

\newcommand{\degree}{\mbox{\it deg}\,}

\newcommand{\itXi}{{\it \Xi}}

\begin{document}

\enlargethispage{23cm}

\begin{titlepage}

$ $

\vspace{-1.5cm} 

\noindent\hspace{-1cm}
\parbox{6cm}{January 2000}\
   \hspace{6.5cm}\
   \parbox{5cm}{math.AG/0002031}

\vspace{2cm}

\centerline{\large\bf
 On the splitting type of an equivariant vector bundle}
\vspace{1ex}
\centerline{\large\bf over a toric manifold}

\vspace{2cm}

\centerline{\large Chien-Hao Liu$\,^1$
              \hspace{1em}and\hspace{1em}
                   Shing-Tung Yau$\,^2$}
\vspace{1.1em}
\centerline{\it Department of Mathematics}
\centerline{\it Harvard University}
\centerline{\it Cambridge, MA 02138}
%

\vspace{3em}

\begin{quotation}
\centerline{\bf Abstract}
\vspace{0.3cm}

\baselineskip 12pt  
{\small
 From the work of Lian, Liu, and Yau on "Mirror Principle", in the
 explicit computation of the Euler data $Q=\{\,Q_0, Q_1,\,\cdots\,\}$
 for an equivariant concavex bundle ${\cal E}$ over a toric manifold,
 there are two places the structure of the bundle comes into play:
  (1) the multiplicative characteric class $Q_0$ of $V$ one starts
      with, and
  (2) the splitting type of ${\cal E}$.
 Equivariant bundles over a toric manifold has been classified
 by Kaneyama, using data related to the linearization of the toric
 action on the base toric manifold. In this article, we relate the
 splitting type of ${\cal E}$ to the classifying data of Kaneyama.
 From these relations, we compute the splitting type of a couple of
 nonsplittable equivariant vector bundles over toric manifolds that
 may be of interest to string theory and mirror symmetry.
 A code in Mathematica that carries out the computation of some of
 these examples is attached.
} 
\end{quotation}

\bigskip

\baselineskip 12pt
{\footnotesize
 } 

\bigskip

\noindent {\small
MSC number 1991: 14M25, 14Q99, 55R91, 14J32, 81T30.
} 

\bigskip

{\footnotesize
\noindent{\bf Acknowledgements.}
 We would like to thank
  Ti-Ming Chiang,
  Shinobu Hosono, Yi Hu,
  Albrecht Klemm,
  Bong H.\ Lian, Kefeng Liu, 
  Jason Starr, 
  and Richard Thomas
 for valuable conversations, discussions, and inspirations
 at various stages of the work.
 C.H.L.\ would like to thank in addition
  Hung-Wen Chang and Ling-Miao Chou
 for discussions of the code.
 The work is supported by DOE grant DE-FG02-88ER25065 and NSF grant
 DMS-9803347.
} 

\noindent
\underline{\hspace{20em}}

$^1${\footnotesize E-mail: chienliu@math.harvard.edu}

$^2${\footnotesize E-mail: yau@math.harvard.edu}

\end{titlepage}

\newpage
$ $

\vspace{-4em}  

\centerline{\sc  Splitting Type of Equivariant Bundles}

\vspace{2em}

\baselineskip 14pt  

\begin{flushleft}
{\Large\bf 0. Introduction and outline.}
\end{flushleft}

\begin{flushleft}
{\bf Introduction.}
\end{flushleft}
From the work of Lian, Liu, and Yau on "Mirror Principle", in the
explicit computation of the Euler data $Q=\{\,Q_0, Q_1,\,\cdots\,\}$
for an equivariant concavex bundle ${\cal E}$ over a toric manifold,
there are two places the structure of the bundle comes into play:
 (1) the multiplicative characteric class $Q_0$ of $V$ one starts
     with, and
 (2) the splitting type of ${\cal E}$.
Equivariant bundles over a toric manifold has been classified by
Kaneyama, using data related to the linearization of the toric
action on the base toric manifold. The purpose of these notes is
to relate the splitting type of ${\cal E}$ to these classification
data of Kaneyama.

In Sec.\ 1, we provide some general backgrounds and notations for
this article. Some basics of equivariant vector bundles over 
toric manifolds are provides in Sec.\ 2. In Sec.\ 3, we discuss how
the splitting type of an equivariant vector bundle over a toric
manifold, if exists, can be obtained from the bundle data.
In Sec.\ 4, we give two classes of examples. In Example 4.1,
we determine which non-decomposable equivariant rank $2$ bundles
over $\CP^2$ admit a splitting type and work out their splitting
type. In Example 4.2, we discuss the splitting type of
tangent/cotangent bundles of toric manifolds. We start with the
splitting type for the (co)tangent bundle of $\CP^n$ and then turn to
the case of toric surfaces. The isomorphism classes of the latter are
coded in a weighted circular graph. From these weights, one can
decide whose (co)tangent bundle admits a splitting type. Since all
toric surfaces arise from consecutive equivariant blowups of either
$\CP^2$ or one of the the Hirzebruch surfaces ${\Bbb F}_a$ and how
the weights on the weighted circular graph behavior under equivariant
blowup is known, the task of deciding which (co)tangent bundle admits
a splitting type and determining them for those that admit one can be
done with the aid of computer. For the interest of string theory and
mirror symmetry, from the toric surfaces arising from equivaraint
blowups of $\CP^2$ up to $9$ points, we sort out those whose
(co)tangent bundle admits a splitting type. The splitting type for
their tangent bundle is also computed and listed. A package in
Mathematica that carries out this computation is attached for
reference.

Overall, this is part of the much bigger ambition of toric
mirror symmetry computation via Euler data, as discussed in
[L-L-Y1, L-L-Y2, L-L-Y3]. We leave the application of the current
article to this goal for another work.

\bigskip

\begin{flushleft}
{\bf Outline.}
\end{flushleft}
{\small
\baselineskip 11pt  
\begin{quote}
 1. Essential backgrounds and notations.

 2. Equivariant vector bundles over a toric manifold and their
    classifications.
      
 3. The splitting type of a toric equivariant bundle.

 4. The splitting type of some examples.

 5. Remarks and issues for further study.

 Appendix. The computer code.
\end{quote}
} 

\newpage
\baselineskip 14pt  

\section{Essential backgrounds and notations for physicists.}

In this section, we collect some basic facts and notations that will
be needed in the discussion. The part that is related to equivariant
vector bundles over toric manifolds is singled out in Sec.\ 2.
Readers are referred to the listed literatures for more details.

\bigskip

\noindent $\bullet$
{\bf Toric geometry.}
([A-G-M], [C-K], [Da], [Ew], [Fu], [Gre], [G-K-Z], [Ke],
and [Od1,Od2].)
Physicists are referred particularly to [A-G-M] or [Gre] for a nice
expository of toric geometry. Let us fix the terminology and
notations here and refer the details to [Fu].

\bigskip

{\it Notation}$\,$:
\begin{quote}
  $N\cong{\Bbb Z}^n$: a lattice;

  $M=\Hom(N,{\Bbb Z})\,$: the dual lattice of $N$;

  $T_N=\Hom(M,{\Bbb C}^{\ast})\,$: the (complex) $n$-torus;

  $\Sigma\,$: a fan in $N_{\Bbb R}$;

  $X_{\Sigma}\,$: the toric variety associeted to $\Sigma$;

  $\Sigma(i)\,$: the $i$-skeleton of $\Sigma$;

  $U_{\sigma}\,$: the local affine chart of $X_{\Sigma}$ associated
     to $\sigma$ in $\Sigma$;

  $x_{\sigma}\in U_{\sigma}\,$: the distinguished points associated
     to $\sigma$;

  $O_{\sigma}\,$: the $T_N$-orbit of $x_{\sigma}$ under the
     $T_N$-action on $X_{\Sigma}$;

  $V(\sigma)\,$: the orbit closure of $O_{\sigma}$;

  $M(\sigma)=\sigma^{\perp}\cap M$.
\end{quote}

\bigskip

Recall that points $v$ in the interior of $\sigma\cap N$ represent
one-parameter subgroups $\lambda_v$ in ${\Bbb T}_N$ such that
$\lim_{z\rightarrow 0}\lambda_v(z)=x_{\tau}$. Recall also that 
the normal cones associated to a polyhedron $\Delta$ with vertices
in $M$ form a fan in $N_{\scriptsizeBbb R}$, called the normal 
fan of $\Delta$. This determines a projective toric variety.

\bigskip

\noindent $\bullet$
{\bf Toric surface.} ([Od2].)
Any complete nonsingular toric surface is obtained from consecutive
equivariant blowups of either $\CP^2$ or one of the Hurzebruch
surfaces ${\Bbb F}_a$ at $T_N$ fixed points. Indeed, one has a
complete classification of them as follows:

\bigskip

\noindent
{\bf Fact 1.1 [toric surface].} ([Od2].) {\it
 The set of isomorphism classes of complete nonsingular toric
 surfaces $X_{\Sigma}$ is in one-to-one corresponce with the set of
 equivalent classes of weighted circular graphs
 $w=(w_1, \,\cdots,\,w_s)$ (under rotation and reflection)
 of the following form: {\rm (\sc Figure 1-1.)}
 \begin{quote}
  \hspace{-1.9em}(1)
  The circular graph having $3$ vertices with weights $1,\,1,\,1$.

  \hspace{-1.9em}(2)
  The circular graph having $4$ vertices with weights in circular
  order $0,\,a,\,0,\,-a$.

  \hspace{-1.9em}(3)
  The weighted circular graphs with $s\ge 5$ vertices that is
  obtained from one with $(s-1)$ vertices by adding a vertex of
  weight $1$ and reducing the weight of each of its two adjacent
  vertices by $1$.
 \end{quote}
} 
 \begin{figure}[htbp]
  \setcaption{{\sc Figure 1-1.}
   \baselineskip 14pt
    The weighted circular graphs that labels the isomorphism classes
    of complete nonsingular toric surfaces.
  } 
  \centerline{\psfig{figure=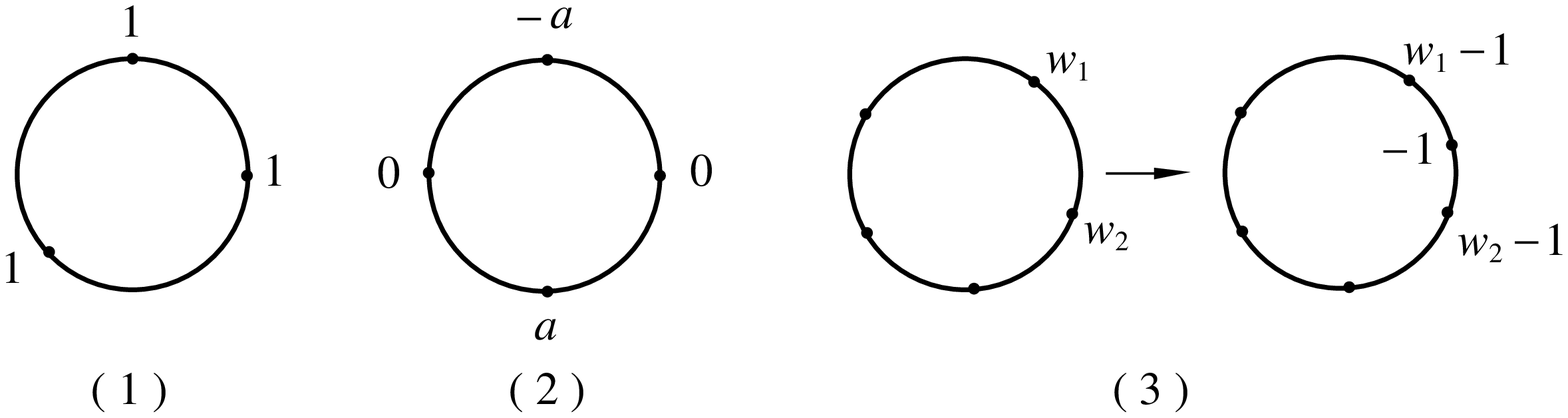,width=13cm,caption=}}
 \end{figure}
Let $\Sigma(1)=(v_1,\,\cdots,\,v_s)$ in, say, counterclockwise order
in $N_{\scriptsizeBbb R}$, then for each $i$, there exists a unique
integer $a_i$ such that $v_{i-1}+v_{i+1}+a_iv_i=0$
(here, $s+1\equiv 1$).
The correspondence is then given by
$X_{\Sigma}\,\mapsto\, (a_1,\,\cdots,\, a_s)$.

\bigskip

\noindent $\bullet$
{\bf Line bundles with positive/negative $c_1$.}
([C-K], [Hi], and [Re].)
Recall that the K\"{a}hler cone in $H^2(X_{\Sigma}, {\Bbb R})$
consists of all the K\"{a}hler classes of $X_{\Sigma}$ and the
Mori cone of $X_{\Sigma}$ consists of all the classes in
$H_2(X_{\Sigma}, {\Bbb R})$ representable by effective $2$-cycles.
From [Re], the Mori cone of $X_{\Sigma}$ is generated by $V(\tau)$,
$\tau\in \Sigma(n-1)$. We call a class
$\omega\in H^2(X_{\Sigma}, {\Bbb R})$ {\it positive}
(resp.\ {\it negative}), in notation $\omega>0$ (resp.\ $\omega<0$),
if $\omega$ (resp.\ $-\omega$) lies in the K\"{a}hler cone of
$X_{\Sigma}$. The fact that the K\"{a}hler cone and the Mori cone
of a complete toric manifold are dual to each other gives us 
a criterion for a line bundle over $X_{\Sigma}$ to have positive
$c_1$:

\bigskip

\noindent
{\bf Fact 1.2 [positive/negative line bundle].} {\it
 Given an $n$-dimensional toric manifold $X_{\Sigma}$.
 Let $L\in\,\Pic(X_{\Sigma})$ be a line bundle over $X_{\Sigma}$ and
 $D$ be the associated divisor class as an element in
 $H_{2n-2}(X_{\Sigma},{\Bbb Z})$. Then $c_1(L)>0$ (resp.\ $<0$) if
 and only if $D\cdot V(\tau)>0$ (resp.\ $<0$) for all
 $\tau\in \Sigma(n-1)$.
} 

\bigskip

\noindent $\bullet$
{\bf The augmented intersection matrix.}
([Fu].)
Given an $n$-dimensional complete nonsingular fan $\Sigma$. Let
$\Sigma(1)=\{\, v_1,\,\cdots,\, v_J\,\}$ and
$\Sigma(n-1)=\{\, \tau_1, \, \cdots,\, \tau_I \,\}$.
Let $A_1$ and $A_{n-1}$ be respectively the first and $(n-1)$-th
Chow group of $X_{\Sigma}$. Then $A_1$ is generated by
$V(\tau_i)$, $i=1,\,\cdots,\, I$, and $A_{n-1}$ is generated by
$D(v_j)$, $j=1,\,\cdots, \, J$. There is a nondegenerate pairing
$A_1\times A_{n-1}\rightarrow {\Bbb Z}$ by taking the intersection
number. Let $Q$ be the $I\times J$ matrix whose $(i,j)$-entry is
the intersection number $V(\tau_i)\cdot D(v_j)$. Since the generators
for $A_1$ and $A_{n-1}$ used here may not be linearly independent,
we shall call $Q$ the {\it augmented intersection matrix}
(with respect to the generators). Explicitly, $Q$ can be determined
as follows.

Let $\tau_i=[v_{j_1},\,\cdots,\,v_{j_{n-1}}]\in\Sigma(n-1)$. Then
$\tau_i$ is the intersection of two $n$-cones
$$
 \sigma_1=[v_{j_1},\,\cdots,\,v_{j_{n-1}}, v_{j_n}]
 \hspace{1em}\mbox{and}\hspace{1em}
 \sigma_2=[v_{j_1},\,\cdots,\,v_{j_{n-1}}, v_{j^{\prime}_n}]
$$
in $\Sigma$. These vertices in $\sigma_1\cup\sigma_2$ satisfy
a linear equation of the form
$$
 v_{j_n}\, +\, v_{j^{\prime}_n}\, +\, a_1v_{j_1}\, 
   +\, \cdots\, +\, a_{n-1}v_{j_{n-1}}\;
  =\; 0\,,
$$
for some unique integers $a_1,\,\cdots,\,a_n$ determined by
$\sigma_1\cup\sigma_2$  . In terms of this, the $i$-th row of $Q$
is simply the coefficient (row) vector of the above equation.
I.e.\ $V(\tau_i)\cdot D(v_{j_k})=a_k$ for $k=1,\,\cdots,\, n-1$;
$V(\tau_i)\cdot D(v_{j_n})=V(\tau_i)\cdot D(v_{j^{\prime}_n})=1$;
and $V(\tau_i)\cdot D(v_j)=0$ for all other $j$.

\bigskip

\noindent $\bullet$
{\bf Cox homogeneous coordinate ring of a toric manifold.}
([Co] and [C-K], also [Au] and [Do].)
Let $\Sigma$ be a fan in ${\Bbb R}^n$ with $\Sigma(1)$ generated
by $\{v_1,\,\cdots,\,v_a\}$ and $A_{n-1}(X_{\Sigma})$ be the Chow
group of $X_{\Sigma}$. Let $(z_1,\,\cdots,\,z_a)$ be the coordinates
of ${\Bbb C}^a$. For $\sigma=[v_{j_1},\,\cdots,\,v_{j_k}]\in\Sigma$,
denote by $z^{\widehat{\sigma}}$ the monomial from
$(z_1\,\cdots\,z_a)/(z_{j_1}\,\cdots\,z_{j_k})$ after cancellation.
Then $X_{\Sigma}$ can be realized as the geometric quotient
$$
 X_{\Sigma}\;
 =\;({\Bbb C}^{\Sigma(1)}-Z(\Sigma))/\mbox{\raisebox{-.4ex}{$G$}}\,,
$$
where $Z(\Sigma)$ is the exceptional subset 
$\{\,(z_1,\cdots,z_a)\,|\, z^{\widehat{\sigma}}=0
             \hspace{1ex}\mbox{for all $\sigma$ in $\Sigma$}\,\}$
in ${\Bbb C}^a$ and $G$ is the group
$\Hom_{\scriptsizeBbb Z}(A_{n-1}(X_{\Sigma}), {\Bbb C}^{\ast})$
that acts on ${\Bbb C}^a$ via the embedding in $({\Bbb C}^{\ast})^a$,
obtained by taking $\Hom(\,\cdot\,,{\Bbb C}^{\ast})$ of the
following exact sequence
$$
 \begin{array}{cccccccccl}
  0 & \longrightarrow  & M  & \longrightarrow  & {\Bbb Z}^a
    & \longrightarrow  & A_{n-1}(X_{\Sigma})
                              & \longrightarrow  & 0  & \\[.6ex]
    & &  m  & \longmapsto  & (m(v_1), \,\cdots,\,m(v_a))
  & & & & &. \\[.6ex]
 \end{array}
$$
More facts will be recalled along the way when we need them.
Their details can be found in Sec.\ 1-3 in [Co].

\bigskip

\section{Basics of equivariant vector bundles over toric manifolds.}

Equivariant vector bundles over a toric manifold have been classified
by Kaneyama and Klyachko independently, using differently sets of
data ([Ka1] and [Kl]). In this article, we use Kaneyama's data as the
starting point to compute splitting types. Some necessary facts
from [Ka1] and [Kl] are summarized below with possibly slight
modification/rephrasing to make the geometric picture more
transparent.

\bigskip

\begin{flushleft}
{\bf Equivariant vector bundles over a toric manifold.}
\end{flushleft}
A vector bundle ${\cal E}$ over a toric manifold $X_{\Sigma}$ is
equivariant if $g^{\ast}{\cal E}={\cal E}$ for all $g\in {\Bbb T}_N$.
An equivariant bundle is linearizable if the action on the base
can be lifted to a fiberwise linear action on the total space of
the bundle.

\bigskip

\noindent
{\bf Fact 2.1 [linearizability].} {\it
 Every equivaraint vector bundle ${\cal E}$ over a toric manifold
 $X_{\Sigma}$ is linearizable.
} 

\bigskip

In general, a bundle can be described by its local trivializations
and the pasting maps. For an equivariant vector bundle ${\cal E}$
over a toric manifold $X_{\Sigma}$, these data can be integrated
with the linearization of the toric action.
\begin{quote}
 \hspace{-1.9em}(1)
 {\it Local trivializations}$\,$: \hspace{1ex}
 Over each invariant affine chart $U_\sigma$ for $\sigma\in\Sigma$,
 the bundle splits:
 $$
  {\cal E}|_{U_{\sigma}}\;
   =\;\oplus_{\chi}\,U_{\sigma}\times E_{\sigma}^{\chi}\,,
 $$
 where $E_{\sigma}^{\chi}$ is the representation of $T_N$ associated
 to the weight $\chi\in M$.

 \hspace{-1.9em}(2)
 {\it The pasting maps}$\,$: \hspace{1ex}
 Over each orbit
 $O_{\tau}\hookrightarrow U_{\sigma_1}\cap U_{\sigma_2}$,
 the pasting map
 $\varphi_{\sigma_2\sigma_1}:
        U_{\sigma_1}\cap U_{\sigma_2}\rightarrow \GL(r,{\Bbb C})$
 is determined by its restriction at a point, say, $x_{\tau}$, in
 $O_{\tau}$, due to the equivariant requirement. The pasting maps
 over different orbits in $U_{\sigma_1}\cap U_{\sigma_2}$ are
 related to each other by the {\it holomorphicity requirement} that
 $\varphi_{\sigma_2\sigma_1}:
   U_{\sigma_1}\cap U_{\sigma_2} \rightarrow \GL(r,{\Bbb C})$
 must be holomorphic for every pair of $\sigma_1$, $\sigma_2$ in
 $\Sigma$. This implies that indeed $\varphi_{\sigma_2\sigma_1}$
 is completely determined by its restriction to a point, say, $x_0$,
 in the dense open orbit $U_0\subset U_{\sigma_1}\cap U_{\sigma_2}$.
 Together with the local splitting propertities in (1) above,
 in fact $\varphi_{\sigma_2\sigma_1}$ is a regular matrix-valued
 function on $U_{\sigma_1}\cap U_{\sigma_2}$. Pasting maps also have
 to satisfy the {\it cocycle condition}$\,$:
 $\varphi_{\sigma_3\sigma_1}
  =\varphi_{\sigma_3\sigma_2}\circ \varphi_{\sigma_2\sigma_1}$
 over $U_{\sigma_1}\cap U_{\sigma_2}\cap U_{\sigma_3}$
 for every triple of $\sigma_1$, $\sigma_2$, $\sigma_3$ in $\Sigma$.
 This implies that the full set of pasting maps between affine
 charts of $X_{\Sigma}$ is determined by the set of pasting maps
 $\varphi_{\sigma_2\sigma_1}$ with
 $\sigma_1,\,\sigma_2\in\Sigma(n)$ that satisfy the cocycle
 condition.
\end{quote}
These observations lead to Kaneyama's data for equivariant vector
bundles over $X_{\Sigma}$.

\bigskip

\begin{flushleft}
{\bf The bundle data and the classification after Kaneyama.}
\end{flushleft}

\noindent
(a) {\it The data of local trivialization}$\,$:
    {\it a collection of weight systems}.

\medskip

For $\sigma\in\Sigma(n)$, $x_{\sigma}$ is a fixed point of the toric
action; thus ${\cal E}_{x_{\sigma}}$ is an invariant fiber of the
lifted toric action. Associated to the representation of ${\Bbb T}^n$
on ${\cal E}_{x_{\sigma}}$ is the weight system
${\cal W}_{\sigma}\subset M$. ${\cal W}_{\sigma}$ determones the
local trivialization of ${\cal E}|_{U_{\sigma}}$:
${\cal E}|_{U_{\sigma}}=\oplus_{\chi\in{\cal W}_{\sigma}}
                               U_{\sigma}\times E^{\chi}_{\sigma}$.

\bigskip

\noindent
(b) {\it The data of pasting}$\,$:
    {\it net of weight systems and pasting maps}.
     
\medskip

These weight systems must satisfy a compatibility condition as
follows. Let $\tau=\sigma_1\cap\sigma_2\in\Sigma(n-1)$ be the common
codimension-$1$ wall of two maximal cones $\sigma_1$ and $\sigma_2$,
then the stabilizer $\Stab(x_{\tau})$ of $x_{\tau}$ is an
$(n-1)$-subtorus in ${\Bbb T}_N$ associated to the sublattice in $N$
spanned by $\tau$. Associated to the representation of
$\Stab(x_{\tau})$ on the fiber ${\cal E}_{x_{\tau}}$ is a weight
system ${\cal W}_{\tau}\subset M/M(\tau)$, where
$M(\tau)=\tau^{\perp}\cap M$. The projection map
$M\rightarrow M/M(\tau)$ induces the maps
$$
 {\cal W}_{\sigma_1}\;
   \stackrel{\pi_{\tau\sigma_1}}{\longrightarrow} \;
 {\cal W}_{\tau}\;
   \stackrel{\pi_{\tau\sigma_2}}{\longleftarrow}\;
 {\cal W}_{\sigma_2}
$$
between weight systems. Since they correspond to the refinement of
the $\Stab(\tau)$-weight spaces to the ${\Bbb T}_N$-weight spaces,
the holomorphicity requiremnent of equivariant pasting maps implied
that both of these maps are surjective. Thus,
$\{\,{\cal W}_{\sigma}\,|\,\sigma\in\Sigma(n)\,\}$
form a net of weight systems.
{\sc Figure 2-1} indicates how the net of weight systems may look
like for $\Sigma$ that comes from the normal cone of a strong
convex polyhedron $\Delta$ in $M$.
 \begin{figure}[htbp]
  \setcaption{{\sc Figure 2-1.}
   \baselineskip 14pt
    For $\Sigma$ being a normal fan of a convex polyhedron in $M$,
    the weight system ${\cal W}$ can be realized as a collection
    of vectors (weights) at the vertices and the barycenter of the
    edges of $\Delta$. The compatibility condition is translated
    into the condition that, for vertices $v_{\sigma_1}$ and
    $v_{\sigma_2}$ connected by an edge $e_{\tau}$ of $\Delta$,
    the three sets ${\cal W}_{\sigma_1}$, ${\cal W}_{\tau}$, and
    ${\cal W}_{\sigma_2}$ have to match up under the projection
    along the direction parallel to $e_{\tau}$.
  } 
  \centerline{\psfig{figure=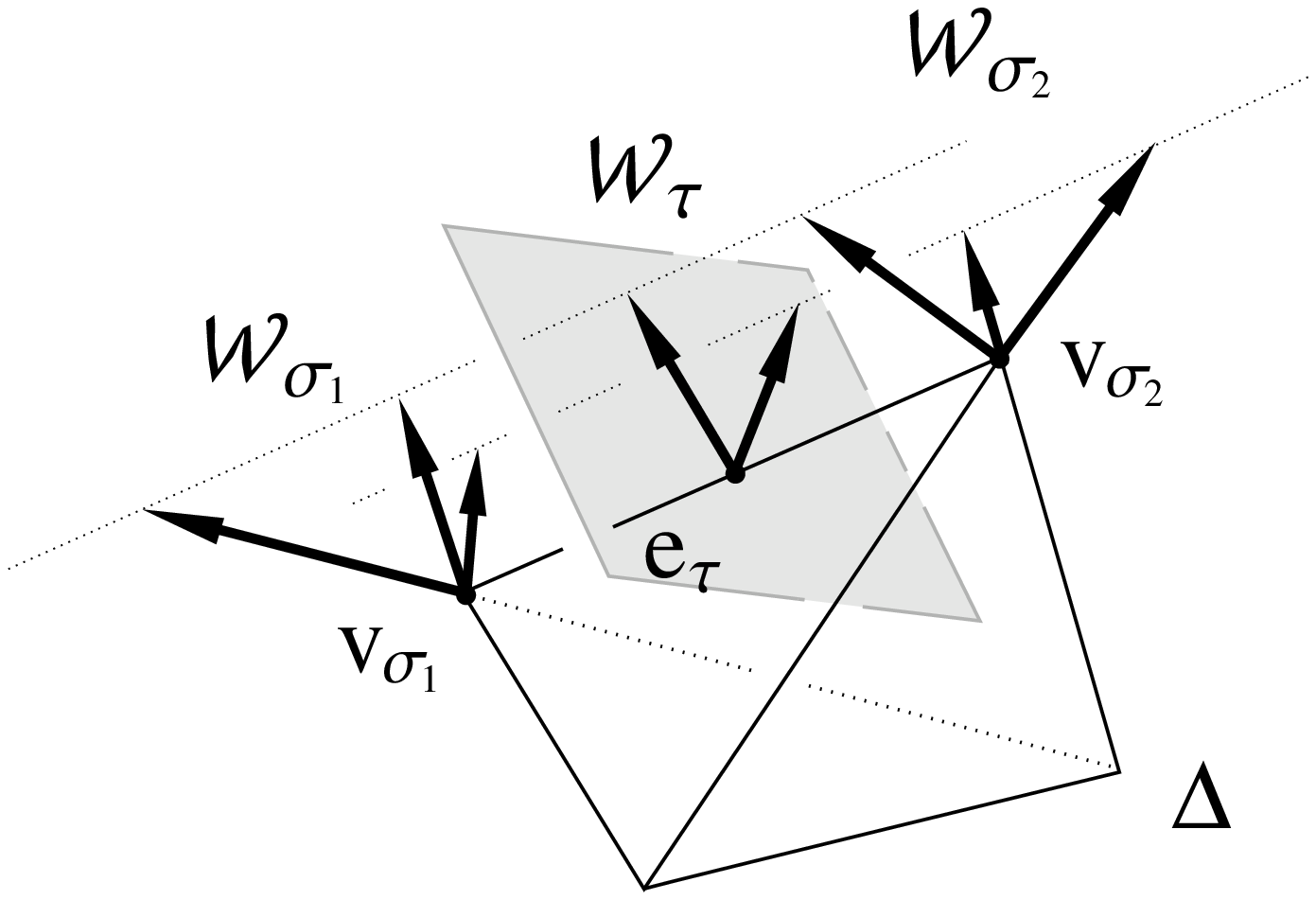,width=13cm,caption=}}
 \end{figure}

The equivariant pasting maps are given by a map 
$$
 P\; :\; \Sigma(n)\times\Sigma(n)\; \longrightarrow\;
            \GL(r,{\Bbb C})
$$
that satisfies
$$
 P(\sigma_3,\sigma_2)\,P(\sigma_2,\sigma_1)\;
 =\; P(\sigma_3,\sigma_1)\,.
$$
$P$ gives a set of compatible pasting maps for the fiber
${\cal E}_{x_0}$ in different ${\cal E}|_{U_{\sigma}}$ with respect
to the bases given by the weight space decomposition. Holomorphicity
condition for its equivariant extension to over
$U_{\sigma_1}\cap U_{\sigma_2}$ requires that:
\begin{quote}
 For any $\tau=\sigma_1\cap\sigma_2\in\Sigma(n-1)$, let
 ${\cal W}_{\sigma_i}
  =(\chi_{\sigma_i 1},\,\cdots,\,\chi_{\sigma_i r})$,
 written with multiplicity given by the dimension of the
 corresponding weight space. Then $P(\sigma_2,\sigma_1)_{ij}=0$
 if $\chi_{\sigma_2 i}-\chi_{\sigma_1 j}\in M-\tau^{\vee}$.
\end{quote}

\bigskip

\noindent
(c) {\it Equivalence of the bundle data.}
 
\medskip

Given $\Sigma$, two bundle data $({\cal W}, P )$ and
$({\cal W}^{\prime}, P^{\prime})$ are said to be equivariant
if ${\cal W}={\cal W}^{\prime}$ and there is a map
$\rho:\Sigma(n)\rightarrow\GL(r,{\Bbb C})$ such that
$P^{\prime}(\sigma_2,\sigma_1)=
           \rho(\sigma_2)P(\sigma_2,\sigma_1) \rho(\sigma_1)^{-1}$.
Equivariant data determine isomorphic linearized equivariant
vector bundles over $X_{\Sigma}$.

\bigskip

\section{The splitting type of an equivariant vector bundle.}

Recall first a theorem of Grothendieck ([Gro] Theorem 2.1), which
says that any holomorphic vector bundle over $\CP^1$ splits into
the direct sum of a unique set of line bundles. The following
definition follows from [L-L-Y2]:

\bigskip

\noindent
{\bf Definition 3.1 [splitting type].} {
 Let ${\cal E}$ be an equivariant vector bundle of rank $r$ over
 a toric manifold $X_{\Sigma}$. Suppose that there exist nontrivial
 equivariant line bundles $L_1,\,\cdots,\, L_r$ over $X_{\Sigma}$
 such that each $c_1(L_i)$ is either $\ge 0$ or $<0$ and that the
 restriction ${\cal E}|_{V(\tau)}$ is isomorphic to the direct sum
 $(\oplus_{i=1}^r\,L_i)|_{V(\tau)}$ for any $\tau\in\Sigma(n-1)$.
 Then $\{L_1,\,\cdots,\, L_r\}$ is called a {\it splitting type} of
 ${\cal E}$.
} 

\bigskip

\noindent
{\bf Definition 3.2 [system of splitting numbers].} {
 For each $\tau\in\Sigma(n-1)$, suppose that
 $$
  {\cal E}|_{V(\tau)}\;
   =\;\oplus_{i=1}^r\,{\cal O}(d^{\tau}_i)
   \hspace{1em}\mbox{with}\hspace{1em}
  d^{\tau}_1\,\ge\, d^{\tau}_2\, \ge\, \cdots\,\ge\, d^{\tau}_r\,.
 $$
 From Grothendieck's theorem, $(d^{\tau}_1,\,\cdots,\,d^{\tau}_r)$
 is uniquely determined by ${\cal E}$. We shall call the set
 $$
  \itXi({\cal E})\;
   =\; \{\, (d^{\tau}_1, \,\cdots\,, d^{\tau}_r)\,|\,
                                           \tau\in\Sigma(n-1) \,\}
 $$
 the {\it system of splitting numbers associated to} ${\cal E}$.
} 

\bigskip

To compute the splitting types of ${\cal E}$, we first extract the
bundle data of ${\cal E}|_{V(\tau)}$ from that of ${\cal E}$ and
then compute $\itXi({\cal E})$ from the bundle data of
${\cal E}|_{V(\tau)}$ by weight bootstrapping. Using these numbers,
one can then determine all the splitting types of ${\cal E}$ by the
augmented intersection matrix $Q$ associated to $\Sigma$.
Let us now turn to the details.

We shall assume that the rank $r$ of ${\cal E}$ $\ge 2$.

\bigskip

\begin{flushleft}
{\bf The bundle data of ${\cal E}|_{V(\tau)}$.}
\end{flushleft}
Let $\tau=\sigma_1\cap\sigma_2\in\Sigma(n-1)$,
$E={\cal E}|_{V(\tau)}$, $U_1=U_{\sigma_1}\cap V(\tau)$, and
$U_2=U_{\sigma_1}\cap V(\tau)$. Let $v_{\tau}$ be a lattice point
in the interior of $\tau$, then the pasting map
$\varphi_{21}(x_{\tau})$ for $E$ from $E|_{U_1}$ to $E|_{U_2}$
over $x_{\tau}$ is given by
$\lim_{z\rightarrow 0} (\lambda_{v_{\tau}}(z)
                P(\sigma_2,\sigma_1)\lambda_{v_{\tau}}(z)^{-1})$,
where $z\in{\Bbb C}^{\ast}$. Since $\Stab(V(\tau))$ acts on $E$
via the $T_N$-action on $E$ and the pasting map for $E$
is $T_N$-equivariant and, hence, commutes with the
$\Stab(V(\tau))$-action, $E$ can be decomposed into
a direct sum of $\Stab(V(\tau))$-weight subbundles:
$$
 E\;=\;\oplus_{\chi\in{\cal W}_\tau}\, E^{\chi}\,.
$$
Indeed, for $\chi\in{\cal W}_{\tau}$,
$E^{\chi}|_{U_1}=\oplus_{\chi^{\prime}\in
                \pi_{\tau\sigma_1}^{-1}(\chi)}\,
                         U_1\times E_{\sigma_1}^{\chi^{\prime}}$;
and similarly for $E^{\chi}|_{U_2}$. Thus, up to a permutation
of elements in the basis, we may assume that
$\varphi_{21}(x_{\tau})$ is in a block diagonal form with each
block labelled by a distinct $\chi\in{\cal W}_{\tau}$.
Let us now turn to the weight system for $E$ at $x_{\sigma_1}$
and $x_{\sigma_2}$.

Let $\tau^{\perp}_{\sigma_1}$ be the primitive lattice point of
$\sigma_1^{\vee}\cap \tau^{\perp}$ in $M$ and $v_{\sigma_1}$ be
a lattice point in $N$ such that
$\langle\tau^{\perp}_{\sigma_1}, v_{\sigma_1}\rangle=1$.
Let $\lambda_{v_{\sigma_1}}$ be the corresponding one-parameter
subgroup in ${\Bbb T}_N$. Then $\lambda_{v_{\sigma_1}}$ acts on
$O_{\tau}$ freely and transitively with
$\lim_{z\rightarrow 0}\lambda_{v_{\sigma_1}}(z)\cdot x_{\tau}
                                                  =x_{\sigma_1}$.
Let
$$
 {\cal W}_{\sigma_1}\;
   =\; \{\,\chi_{\sigma_1 1},\,\chi_{\sigma_1 2},\,\cdots\,\}
 \hspace{1em}\mbox{and}\hspace{1em}
 {\cal W}_{\sigma_2}\;
   =\; \{\,\chi_{\sigma_2 1},\,\chi_{\sigma_2 2},\,\cdots\,\}
$$
be the set of ${\Bbb T}_N$-weights at $x_{\sigma_1}$ and
$x_{\sigma_2}$ respectively. Then, as a
$\lambda_{v_{\sigma_1}}$-equivariant bundle with the induced
linearization from the linearization of ${\cal E}$, the
corresponding $\lambda_{\sigma_1}$-weights of $E_{x_{\sigma_1}}$
and $E_{x_{\sigma_2}}$
are given respectively by
$$
 {\cal W}^{{\scriptsizeBbb T}^1}_{\sigma_1}\;
 =\; \{\,\langle\chi_{\sigma_1 1}, v_{\sigma_1}\rangle,\,
  \langle\chi_{\sigma_1 2}, v_{\sigma_1}\rangle,\,\cdots\,\}
 \hspace{1em}\mbox{and}\hspace{1em}
 {\cal W}^{{\scriptsizeBbb T}^1}_{\sigma_2}\;
 =\; \{\,\langle\chi_{\sigma_2 1}, v_{\sigma_1}\rangle,\,
  \langle\chi_{\sigma_2 2}, v_{\sigma_1}\rangle,\,\cdots\,\}\,.
$$
Notice that both sets depend on the choice of $v_{\sigma_1}$;
however, different choices of $v_{\sigma_1}$ will lead only to
an overall shift of
${\cal W}_{\sigma_1}^{{\scriptsizeBbb T}^1}\cup
                   {\cal W}_{\sigma_2}^{{\scriptsizeBbb T}^1}$
by an integer.

\bigskip

\begin{flushleft}
{\bf From bundle data to splitting numbers$\,$:
     weight bootstrapping.}
\end{flushleft}
Recall first the following fact by Grothendieck [Gro]:

\bigskip

\noindent
{\bf Fact 3.3 [Grothendieck].} {\it
 Given a holomorphic vector bundle $E$ of rank $r$ over $\CP^1$. Let
 $E_0=\{0\}\subset E_1\subset\,\cdots\,\subset E_r=E$ be a filtration
 of $E$ such that $E_i/E_{i-1}$ is a line bundle for $1\le i\le r$
 and that the degree $d_i$ of $E_i/E_{i-1}$ form a non-increasing
 sequence. Then $E$ is isomorphic to the direct sum
 $\oplus_{i=1}^r\,(E_i/E_{i-1})$.
} 

\bigskip

\noindent
Following previous discussions and notations, we only need to work
out the splitting numbers for each $E^{\chi}$. Thus, without loss
of generality, we may assume that ${\cal W}_{\tau}={\chi}$ in the
following discussion.

Fix a $v_{\sigma_1}$. Note that in terms of the one-parameter
subgroup $\lambda_{v_{\sigma_1}}$ acting on $V(\tau)$, $x_{\sigma_1}$
has coordinate $0$ while $x_{\sigma_2}$ has coordinate $\infty$.
Let
$$
 E|_{U_1}=\oplus_{i=1}^a U_1\times E_1^{\chi_{1i}}
 \hspace{1em}\mbox{and}\hspace{1em}
 E|_{U_2}=\oplus_{j=1}^b U_2\times E_2^{\chi_{2j}}
$$
be the induced ${\Bbb T}^1$-weight space decomposition of $E|_{U_1}$
and $E|_{U_2}$ respectively, and $\varphi_{12}(x_{\tau})$ be the
pasting map at $x_{\tau}$ from $(E|_{U_1})|_{x_{\tau}}$ to
$(E|_{U_2})|_{x_{\tau}}$.
We assume that $\chi_{1 1}>\,\cdots\,>\chi_{1 a}$ and
$\chi_{2 1}<\,\cdots\,<\chi_{2 b}$. Our goal is now
to work out a filtration of $E$, using the given bundle data,
that satisfies the property in the above fact. 

Let $v$ be a non-zero vector in the fiber $E_{x_{\tau}}$ over
$x_{\tau}$. Then associated to the ${\Bbb T}^1$-orbit
${\Bbb T}^1\cdot v$ of $v$ is a line bundle ${\cal L}_v$ that
contains ${\Bbb T}^1\cdot v$ as a meromorphic section $s_v$.
Let $v_1^{\chi_{1 i^{\prime}}}$ be the lowest ${\Bbb T}^1$-weight
component of $v$ in $E|_{U_1}$ and $v_2^{\chi_{2 j^{\prime}}}$)
be the highest ${\Bbb T}^1$-weight component of $v$ in
$E|_{U_2}$. Then ${\cal L}_v|_{x_{\sigma_1}}$
(resp.\ ${\cal L}_v|_{x_{\sigma_2}}$) lies in
$E_1^{\chi_1 i^{\prime}}$ (resp.\ $E_2^{\chi_{2 j^{\prime}}}$)
and the meromorphic section $s_v$ is holomorphic over $O_{\tau}$
with a zero at $x_{\sigma_1}$ (resp.\ $x_{\sigma_2}$) of order
$\chi_{1 i^{\prime}}$ (resp.\ $-\chi_{2 j^{\prime}}$).
(Here, a zero of order $-k$ means the same as a pole of order $k$.)
This shows that
$$
 {\cal L}_v\;
    =\; {\cal O}(\chi_{1 i^{\prime}}-\chi_{2 j^{\prime}})\,.
$$
Notice that, from the previous discussion, this degree is
independent of the choices of $v_{\sigma_1}$.
We shall now proceed to construct a ${\Bbb T}^1$-equivariant line
subbundle in $E$ that achieves the maximal degree.

Fix a basis for the ${\Bbb T}^1$-weight spaces in the local
trivialization of $E$, then $\varphi_{21}(x_{\tau})$ is expressed
by a matrix $A$ which admits a weight-block decomposition
$A=[\,A_{\chi_{2 j},\chi_{1 i}}\,]_{ji}$. Consider the chain of
submatrices $B_{kl}$ in $A$ that consists of weight
blocks $A_{\chi_{2 j},\chi_{1 i}}$ with $j=k,\,\cdots,\, b$ and
$i=1,\,\cdots,\, l$, as indicated in {\sc Figure 3-1}.
\begin{figure}[htbp]
 \setcaption{{\sc Figure 3-1.}
  \baselineskip 14pt
    The submatrix $B_{kl}$ of $A$, which consists of weight blocks,
    are indicated by the shaded part.
 } 
 \centerline{\psfig{figure=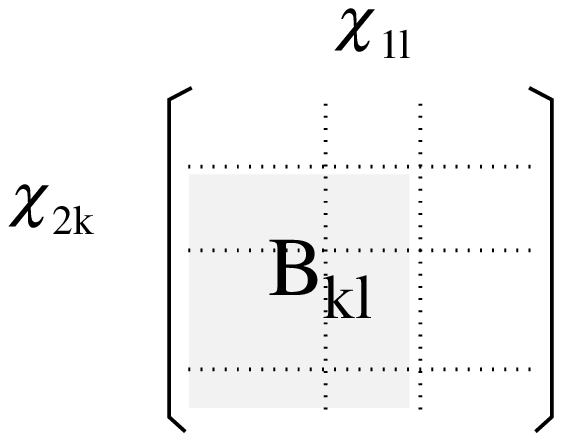,width=13cm,caption=}}
\end{figure}
Each $B_{kl}$ gives the linear map from
$\oplus_{i=1}^l E_1^{\chi_{1 i}}$ to
$\oplus_{j=k}^b E_2^{\chi_{2 j}}$ induced by $A$.
Let $N_{kl}$ be the kernel of $B_{kl}$, as a subspace in
$(E|_{U_1})_{x_{\tau}}$. Define
$$
 {\cal E}_i\;
 =\;E_1^{\chi_{11}}\oplus\,\cdots\,\oplus E_1^{\chi_{1 i}}
 \hspace{1em}\mbox{at $x_{\tau}$}\,,
 \hspace{1em} i=1,\,\cdots,\, a\,.
$$
Then one has the following sequence of filtrations:
$$
 \begin{array}{cccccccll}
  \{0\}  & \subset  & {\cal E}_1  & \subset  & \cdots  & \subset
      & {\cal E}_a & \hspace{-2ex}=\;(E|_{U_1})_{x_{\tau}}& \\[.6ex]
  & & \cup  & & \cdots & & \cup & \\[.6ex]
  & & N_{b1} & \subset  & \cdots  & \subset  & N_{ba} & \\[.6ex]
  & & \cup  & & \cdots & & \cup  & \\[.6ex]
  & & \vdots  & & \vdots  & & \vdots  & \\[.6ex]
  & & \cup  & & \cdots & & \cup  & \\[.6ex]
  & & N_{11}  & \subset  & \cdots  & \subset  & N_{1a} & &.
 \end{array}
$$
Let
$$
 {\cal E}_{ij}\; =\; {\cal E}_i-{\cal E}_{i-1}-N_{ji}\,,
  \hspace{1em}\mbox{for}\hspace{1ex}
  i=1,\,\cdots,\,a\,, j=1,\,\cdots,\, b\,;
$$
then
$$
 (E|_{U_1})_{x_{\tau}}\;=\;\cup_{i,j}\, {\cal E}_{ij}\,.
$$
By construction, if ${\cal E}_{ij}$ is non-empty, then, for
any $v$ in ${\cal E}_{ij}$,
$\degree({\cal L}_v)=\chi_{1 i}-\chi_{2 j}$.
Thus we may define the characteristic number $\mu({\cal E}_{ij})$
for ${\cal E}_{ij}$ non-empty to be $\chi_{1 i}+\chi_{2 j}$.
Now let 
$$
 d_1\; =\; \mbox{\it max}\,\{\, \mu({\cal E}_{ij})\,|\,
              \mbox{${\cal E}_{ij}$ non-empty,
               $i=1,\,\cdots,\, a,\, j=1, \,\cdots,\, b$}\,\}
$$
and 
$v\in {\cal E}_{ij}$ for some $(ij)$ that realizes $d_1$; then
by construction $E_1={\cal L}_v$ is a ${\Bbb T}^1$-equivariant
line subbundle of $E$ that achieves the maximal possible degree.

\bigskip

Suppose the basis for the weight spaces in the local trivialization
of $E$ are given by $(e_{11},\,\cdots,\,e_{1r})$ and
$(e_{21},\,\cdots,\,e_{2r})$ respectively.
Let $v=\sum_{i=1}^{k}c_i e_{1 i}$ with $c_k\ne 0$ in
$E|_{U_1}$ and
$v=\sum_{i=1}^{k^{\prime}}c^{\prime}_i e_{2i}$ with
$c^{\prime}_{k^{\prime}}\ne 0$ in $E|_{U_2}$. Then, by
replacing $e_{1k}$ and $e_{2k^{\prime}}$ by $v$ and
putting it as the first element in the bases, one shows

\bigskip

\noindent
{\bf Lemma 3.4.} {\it
 There exist a new weight space decomposition of $E|_{U_1}$ and
 $E|_{U_2}$ respectively and a choice of the basis for the
 new weight spaces such that $E_1|_{x_{\tau}}$ is spanned by the
 first element in the basis.
} 

\bigskip

\noindent
This renders the pasting map into the form:
$$
 A\; =\; \left[\, \begin{array}{cc}
                    1 & \ast \\[.6ex]
                    0 & A_1
                  \end{array}  \,\right]\,,
$$
where $0$ is the $(r-1)$-dimensional zero vector and $A_1$ is a
nondegenerate $(r-1)\times(r-1)$ matrix.
With respect to the new trivialization, the ${\Bbb T}^1$-weight
spaces and their basis descends then to the quotient $E/E_1$
with pasting map given by $A_1$.

Repeating the discussion $r$ times, one obtains 
${\Bbb T}^1$-equivariant line subbundles
$$
 E_1\subset E\,,\hspace{2ex} E_2/E_1\subset E/E_1\,, \hspace{2ex}
   \cdots\,, \hspace{2ex} E_{r-1}/E_{r-2}\subset E/E_{r-2}\,,
   \hspace{2ex} E_r/E_{r-1}\,.
$$
of maximal degree in each pair and the associated filtration
$$
 E_0=\{0\}\subset E_1 \subset \,\cdots\, \subset E_r=E\,.
$$

\bigskip

\noindent
{\bf Lemma 3.5.} {\it
 The filtration of $E$ obtained above satisfies the condition of
 Grothendieck in Fact 3.3.
} 

\bigskip

\noindent
{\it Proof.}
Since $E_i/E_{i-1}$ is a ${\Bbb T}^1$-equivariant line subbundle
of maximal degree in $E/E_{i-1}$, it must be so also in
$E_{i+1}/E_{i-1}$. Together with the fact that 
$E_{i+1}/E_i=(E_{i+1}/E_{i-1})/(E_i/E_{i-1})$,
we only need to justify the claim for the rank $2$ case.
Thus, assume that $E$ is of rank $2$ and the filtration is given
by $\{0\}\subset E_1\subset E$. By Lemma 3.4, we may choose a basis
compatibe with the ${\Bbb T}^1$-weight space decomposition and
the filtration such that the pasting map $\varphi_{\infty 0}$ and
its inverse are given respectively by the following matrices with
the ${\Bbb T}^1$-weight indicated:
$$
 \begin{array}{c}
   \chi_{1 1}\hspace{1ex}\chi_{1 2}\\[.6ex]
   A\; =\; \left[\, \begin{array}{cc}
                    1 & \ast  \\
                    0 & 1 
                  \end{array}  \,\right]\;
           \begin{array}{l}
             \chi_{2 1} \\ \chi_{2 2}
           \end{array}
 \end{array}
 \hspace{2em}\mbox{\raisebox{-1.7ex}{and}}\hspace{2em}
 \begin{array}{c}
   \hspace{1.5ex}\chi_{2 1}\hspace{1ex}\chi_{2 2}\\[.6ex]
   A^{-1}\; =\; \left[\, \begin{array}{cr}
                    1 & -\ast  \\
                    0 & 1
                  \end{array}  \,\right]\;
           \begin{array}{l}
             \chi_{1 1} \\ \chi_{1 2}
           \end{array}
 \end{array}\,\mbox{\raisebox{-1.7ex}{.}}
$$
This implies that $\degree E_1=\chi_{11}-\chi_{21}$ and
$\degree(E/E_1)=\chi_{12}-\chi_{22}$. The linear independency
of the line bundles associated to the column vectors of $A$ at $z=0$
requires that $\chi_{22}\ge\chi_{21}$. Similarly, the
linear independency of the line bundles associated to the column
vectors of $A^{-1}$ at $z=\infty$ requires that
$\chi_{12}\le\chi_{11}$. Consequently,
$\degree E_1\ge\degree(E/E_1)$. This concludes the proof.

\noindent\hspace{14cm} $\Box$

\bigskip

Consequently, by Fact 3.3 one obtains the decomposition of $E$ as
the direct sum of line bundles and the splitting numbers.

\bigskip

\noindent
{\it Remark 3.6} [$\,${\it weight matching}$\,$].
For $\tau=\sigma_1\cap\sigma_2\in\Sigma(n-1)$, if all the
$\Stab(x_{\tau})$-weight spaces are $1$-dimensional, then up to
a permutation of the elements in bases, the pasting map
$\varphi_{21}(x_{\tau})$ becomes diagonal and the correpondences
${\cal W}_{\sigma_1}\rightarrow {\cal W}_{\tau}
                                \leftarrow {\cal W}_{\sigma_2}$
are bijective. Let
$\chi_{\sigma_1i}\rightarrow\chi_{\tau i}
                                  \leftarrow\chi_{\sigma_2 i}\,$,
$i=1,\,\cdots,\,r$, be the correpondences of weights.
Then, up to a permutation, the splitting number of ${\cal E}$ over
$V(\tau)$ is given by
$$
 \left(\, \langle \chi_{\sigma_11}-\chi_{\sigma_21}\,,\,
                       v_{\sigma_1}\rangle\,,\,\cdots\,,\,
          \langle \chi_{\sigma_1r}-\chi_{\sigma_2r}\,,\,
                       v_{\sigma_1}\rangle
                             \,\right)\,,
$$
where recall that
$\tau^{\perp}_{\sigma_1}=\tau^{\perp}\cap \sigma_1^{\vee}$ and
$\langle\tau^{\perp}_{\sigma_1},v_{\sigma_1}\rangle=1$.

\bigskip

\begin{flushleft}
{\bf From the system of splitting numbers to the splitting types
     of ${\cal E}$ .}
\end{flushleft}
Following the notation in Sec.\ 1, let
$\Sigma(n-1)=\{\,\tau_1,\,\cdots,\,\tau_I\,\}$ and
$$
 \itXi({\cal E})\;
  =\; \{\, (d^{\tau_i}_1, \,\cdots,\, d^{\tau_i}_r )\,|\,
                                       i=1,\,\cdots,\, r\,\}\,.
$$
be the system of splitting numbers associated to ${\cal E}$.
Let $R$ be the $I\times r$ matrix whose $i$-th row is
$(d^{\tau_i}_1, \,\cdots,\, d^{\tau_i}_r )$. Recall the augmented
intersection matrix $Q$ from Sec.\ 1.
Then the problem of finding splitting types of ${\cal E}$ is
equivalent to finding out matrices $R^{\prime}$ obtained
by row-wise permutations of $R$ such that:
\begin{quote}
 \hspace{-1.9em}(1)\hspace{1ex}
 Each column of $R^{\prime}$ has only all positive, all zero, or all
 negative entries.

 \hspace{-1.9em}(2)\hspace{1ex}
 The following matrix linear equation has an integral solution:
 $$
  Q\,X\;=\; R^{\prime}\,,
 $$
 where $X$ is an $J\times r$ matrix.
\end{quote}

Let $X_{kl}$ be the $(k,l)$-entry of $X$. Then associated to the
$r$-many column vectors of the solution matrix $X$ are the line
bundles $L_l$ represented by $\sum_{k=1}^J\,X_{kl}D(v_k)$, for
$l=1,\,\cdots,\, r$. From Fact 1.2 in Sec.1, Condition (1) above for
$R^{\prime}$ means that $c_1(L_l)$ is either $\ge 0$ or $<0$. Such
set of line bundles gives then a splitting type of ${\cal E}$ by
construction. If there exist no such $(R^{\prime}, X)$, then
${\cal E}$ does not admit a splitting type. Finding all such
$R^{\prime}$ and solving the matrix $X$ can be achieved by using
a computer.

\bigskip

\section{The splitting type of some examples. }

In this section, we compute the splitting types of some equivaraint
vector bundles over toric manifolds to illustrate the ideas in
previous sections and also for future use. The details of the toric
manifolds used here can be found in [Fu] and [Od2].

\bigskip

\noindent
{\bf Example 4.1 [equivariant vector bundles of rank 2 over $\CP^2$].}
Recall first ([Fu]) the toric data for $\CP^2$, as illustrated in
{\sc Figure 4-1}(a).
Let ${\cal E}$ be an indecomposable equivariant vector bundle of
rank 2 over $\CP^2=\Proj({\Bbb C}[u_0, u_1, u_2])$.
From [Ka1], ${\cal E}$ is isomorphic to
${\cal E}_{a,b,c,n}={\cal E}(a,b,c)\otimes {\cal O}(n)$
or its dual bundle for some positive integers $a,\,b,\,c$ and
integer $n$, where ${\cal E}(a,b,c)$ is the rank 2 bundle defined
by the exact sequence
$$
 \begin{array}{cccccccccl}
  0 & \longrightarrow & {\cal O}_{{\scriptsizeBbb C}{\rm P}^2}
    & \longrightarrow
    & {\cal O}(a)\oplus{\cal O}(b)\oplus{\cal O}(c)
    & \longrightarrow & {\cal E}(a, b, c)
    & \longrightarrow & 0 &                    \\[,6ex]
  & & 1 & \longmapsto & (\,u_0^a,\, u_1^b,\, u_2^c\,)
    & & & & & .
 \end{array}
$$
From the bundle data of ${\cal E}(a,b,c)$ as worked out in [Ka1],
the weight systems for ${\cal E}(a,b,c)$ at the distinguished points
$x_{\sigma_1},\, x_{\sigma_1}$ and $x_{\sigma_3}$
are given respectively by (cf.\ {\sc Figure 4-1}(b))
$$
 W_1\;=\;\{(a,0),\,(0,b)\}\,, \hspace{2em}
 W_2\;=\;\{(-b,b),\,(-c,0)\}\,, \hspace{2em}
 W_3\;=\;\{(a,-a),\,(0,-c)\}\,.
$$
 \begin{figure}[htbp]
  \setcaption{{\sc Figure 4-1.}
   \baselineskip 14pt
   In (a), the fan and its dual cones for $\CP^2$ are illustrated.
   In (b), the weight systems $W_1,\,W_2,\, W_3$ associated to
   ${\cal E}(a,b,c)$, $a,\,b,\,c$ positive integers, are illustrated.
  } 
  \centerline{\psfig{figure=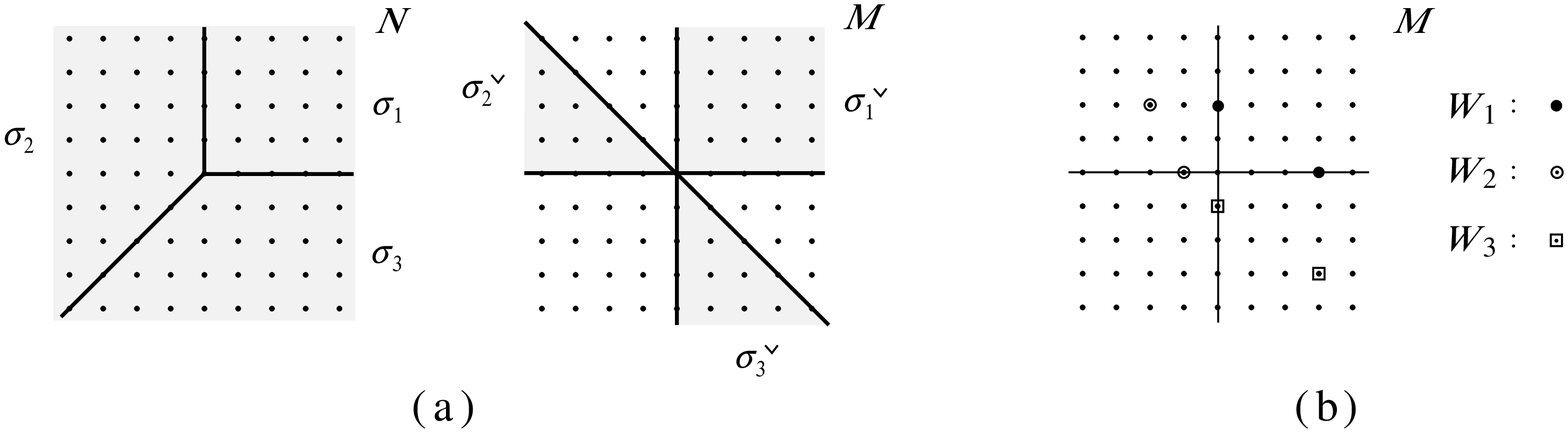,width=13cm,caption=}}
 \end{figure}

Comparing with the toric data for $\CP^2$ and following the
discussions in Sec.\ 3, in particular Remark 3.6, one has
$$
 \begin{array}{l}
 {\cal E}(a,b,c)|_{\overline{x_{\sigma_1}x_{\sigma_2}}}\;
   =\; {\cal O}(a+c)\oplus{\cal O}(b)\,, \hspace{1em}
 {\cal E}(a,b,c)|_{\overline{x_{\sigma_2}x_{\sigma_3}}}\;
   =\; {\cal O}(a+b)\oplus{\cal O}(c)\,, \hspace{1em}
 \mbox{and} \\[.6ex]
 {\cal E}(a,b,c)|_{\overline{x_{\sigma_1}x_{\sigma_3}}}\;
   =\; {\cal O}(b+c)\oplus{\cal O}(a)\,.
\end{array}
$$
Consequently,
$$
 \begin{array}{l}
 {\cal E}|_{\overline{x_{\sigma_1}x_{\sigma_2}}}\;
   =\; {\cal O}(a+c+n)\oplus{\cal O}(b+n)\,, \hspace{1em}
 {\cal E}|_{\overline{x_{\sigma_2}x_{\sigma_3}}}\;
   =\; {\cal O}(a+b+n)\oplus{\cal O}(c+n)\,, \hspace{1em}
 \mbox{and} \\[.6ex]
 {\cal E}|_{\overline{x_{\sigma_1}x_{\sigma_3}}}\;
   =\; {\cal O}(b+c+n)\oplus{\cal O}(a+n)\,.
\end{array}
$$
Thus, up to permutations, the system of splitting numbers associated to
${\cal E}$ is
$$
 \itXi({\cal E})\;
  =\; \{\, (a+c+n, b+n), \, (a+b+n, c+n),\, (b+c+n, a+n)\,\}\,
$$
and
$$
 R\; =\; \left[\, \begin{array}{ccc}
                    a+c+n & b+n \\
                    a+b+n & c+n \\
                    b+c+n & a+n 
                  \end{array} \,\right]\,.
$$
From {\sc Figure} 4-1(a), the augmented intersection matrix $Q$ for
$\CP^2$ is given by
$$
 Q\; =\; \left[\, \begin{array}{rrr}
                    1 & 1 & 1 \\
                    1 & 1 & 1 \\
                    1 & 1 & 1
                  \end{array}   \, \right]\,.
$$
Performing row-wise permutations to $R$, one obtains four possible
$R^{\prime}$, up to overall permutations of column vectors.
Solving the matrix equation $QX=R^{\prime}$, one concludes that

\bigskip

\noindent
{\bf Corollary.} {\it For the indecomposable equivaraint rank $2$
 bundle ${\cal E}_{a, b, c, n}$ over $\CP^2$ to admit a splitting
 type, one must have $a=b=c$. In this case, the splitting type of
 ${\cal E}_{a,a,a,n}$ is unique and is given by
 $(\,{\cal O}(2a+n),\, {\cal O}(a+n)\,)$.
} 

\bigskip

This concludes the example.

\noindent\hspace{14cm} $\Box$

\bigskip

\noindent
{\bf Example 4.2 [(co)tangent bundle of toric manifolds].}
Notice that, since $T_{\ast}X$ and $T^{\ast}X$ are dual to each
other, their splitting types are negative to each other. Thus,
we only need to consider $T_{\ast}X$. Let us compute first the
splitting numbers of $T_{\ast}X_{\Sigma}$.
Let $\tau=\sigma_1\cap\sigma_2\in \Sigma(n-1)$ with
$$
 \sigma_1=[v_1,\,\cdots,\,v_{n-1}, v_n]
 \hspace{1em}\mbox{and}\hspace{1em}
 \sigma_2=[v_1,\,\cdots,\,v_{n-1}, v_{n^{\prime}}]\,.
$$
These vertices in $\sigma_1\cup\sigma_2$ satisfy a linear equation
of the form
$$
 v_{j_n}\, +\, v_{j^{\prime}_n}\, +\, a_1v_{j_1}\,
   +\, \cdots\, +\, a_{n-1}v_{j_{n-1}}\;
  =\; 0\,,
$$
for some unique integers $a_1,\,\cdots,\,a_n$ determined by
$\sigma_1\cup\sigma_2$. Let $(e^1,\,\cdots,\,e^n)$ be the dual
basis in $M$ with respect to $(v_1,\,\cdots,\,v_n)$. Then
$$
 \sigma_1^{\vee}=[e^1,\,\cdots,\,e^{n-1}, e^n]
 \hspace{1em}\mbox{and}\hspace{1em}
 \sigma_2^{\vee}=[e^1-a_1e^n,\,\cdots,\,e^{n-1}-a_{n-1}e^n, -e^n]\,.
$$
Consequently, $\chi_{\sigma_1i}=e^i$ for $i=1,\,\cdots,\,n$,
$\chi_{\sigma_2i}=e^i-a_ie^n$ for $i=1,\,\cdots,\, n-1$; and
$\chi_{\sigma_2 n}=-e^n$. Choosing $v_{\sigma_1}=v_n$ and by the
discussion in Sec.\ 3, one concludes that the splitting number of
$T_{\ast}X_{\Sigma}$ over $V(\tau)$ is given by 
$$
 (\,a_1,\,\cdots,\, a_{n-1},\, 2\,)\,,
$$
up to a permutation. Thus, the system
$\itXi(T_{\ast}X_{\Sigma})$ of spliting numbers associated to
$T_{\ast}X_{\Sigma}$ is already coded in $\Sigma$, as it should be.

In the dual picture, if $X_{\Sigma}$ is projective and, hence,
$\Sigma$ is realized as the normal fan of a strongly convex
polyhedron $\Delta$ in $M$. Let $m_{\sigma}$ be the vertex of
$\Delta$ associated to $\sigma\in\Sigma(n)$. Then ${\cal W}_{\sigma}$
is the set of primitive vectors that generate the tangent cone of
$\Delta$ at $m_{\sigma}$.
If $\tau=\sigma_1\cap\sigma_2\in\Sigma(n-1)$, then
$m_{\sigma_1}$ and $m_{\sigma_2}$ are connected by an edge
$\overline{m_{\sigma_1} m_{\sigma_2}}$ of $\Delta$ that is parallel
to $\tau^{\perp}$. Thus the projection $M\rightarrow M/M(\tau)$ is
given by the projection along the
$\overline{m_{\sigma_1} m_{\sigma_2}}$-direction.
The strong convexity of $\Delta$ implies that
${\cal W}_{\sigma_1}$ and ${\cal W}_{\sigma_2}$ match up
bijectively under this projection. Thus $\itXi(X_{\Sigma})$
can be also read off directly from $\Delta$.

We can now compute the splitting type of some concrete examples.
The result shows that$\,$: {\it Not every tangent bundle of a toric
manifold admits a splitting type}.

\bigskip

\noindent
(a) {\it The projective space $\CP^n$.}
Let $(v_1,\,\cdots,\, v_n)$ be a basis of $N$ and
$v_{n+1}=-(v_1\,\cdots\,v_n)$ and $\Sigma$ be the fan whose maximal
cones are generated by every independent $n$ elements in
$\{\,v_1,\,\cdots\,v_{n+1}\,\}$. Then $\CP^n=X_{\Sigma}$.
By construction, $\Sigma(1)$ is given by
$\{v_1,\,\cdots,\, v_{n+1}\}$ with $v_1+\,\cdots\,v_{n+1}=0$. Thus
the splitting number for $T_{\ast}\CP^n$ over any invariant $\CP^1$
is given by
$$
    (\,2,\,\underbrace{1,\,\cdots,\,1}_{n-1}\,)\,.
$$
The augmented intersection matrix $Q$ has all of its entries equal
to $1$. From this, one concludes that the splitting type of
$\CP^n$ is unique and is given by
$$
 (\,{\cal O}(2),\, \underbrace{{\cal O}(1),\,
                             \cdots,\,{\cal O}(1) }_{n-1}\,)\,.
$$

\bigskip

\noindent
(b) {\it The Hirzebruch surface ${\Bbb F}_a$.}
The toric data for the Hurzebruch ${\Bbb F}_a$ and its weighted
circular graph [Od2] is given in {\sc Figure 4-2}.
 \begin{figure}[htbp]
  \setcaption{{\sc Figure 4-2.}
   \baselineskip 14pt
   The fan for the Hirzebruch surface ${\Bbb F}_a$ and its weights.
  } 
  \centerline{\psfig{figure=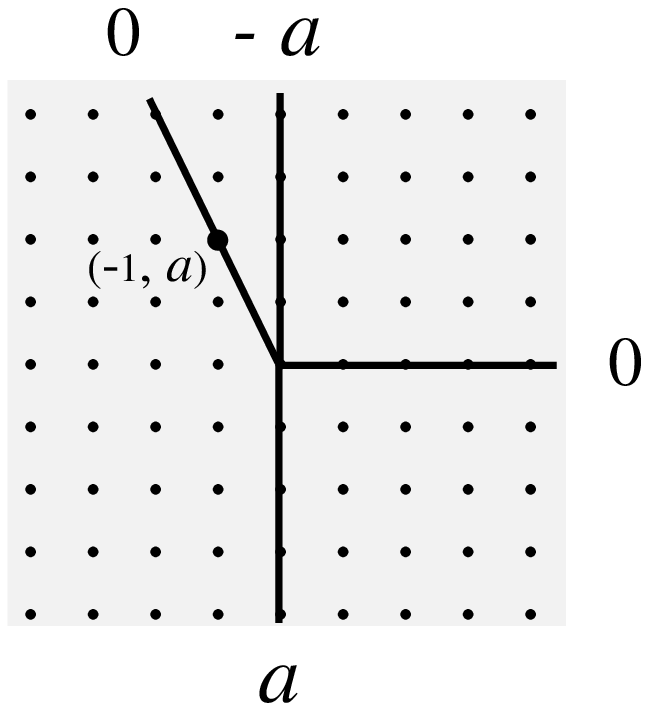,width=13cm,caption=}}
 \end{figure}

Consequently,
$$
 \itXi(T_{\ast}{\Bbb F}_a)\;
  =\; \{\, (2, 0),\,(2, a),\, (2, 0),\, (2, -a)\,\}
$$
and its augmented intersection matrix is given by
$$
 Q\; =\;\left[\, \begin{array}{rrrr}
                   0 & 1 & 0 & 1 \\
                   1 &  a & 1 & 0 \\
                   0 &  1 & 0 & 1 \\
                   1 &  0 & 1 & -a
                 \end{array}   \,\right]\,.
$$
Thus, the only Hirzebruch surface whose tangent bundle can
admit a splitting type is when $a=0$ since the line bundle in a
splitting type must be either positive or negative.
For $a=0$, ${\Bbb F}_0=\CP^1\times\CP^1$. Let
$H^2(\CP^1\times\CP^1, {\Bbb Z})={\Bbb Z}\oplus{\Bbb Z}$ from the
product structure. Then direct computations as in Example 4.1
concludes that $T_{\ast}({\Bbb F}_0)$ admits a unique splitting type
$$
 (\, {\cal O}(2,2),\, {\cal O}({\Bbb F}_0)\,)\,,
$$
where ${\cal O}(2,2)$ is the line bundle associated to $(2,2)$
in $H^2({\Bbb F}_0, {\Bbb Z})$ and ${\cal O}({\Bbb F}_0)$ is the
trivial line bundle.

\bigskip

\noindent
(c) {\it The blowups of $\CP^2$ or ${\Bbb F}_a$.}
Recall from [Fu] and [Od2] that every complete nonsingular toric
surface $X$ is obtained from $\CP^2$ or ${\Bbb F}_a$, $a>0$, by
a succession of blowups at the $T_N$-fixed points. Let
$(a_1,\, \cdots,\, a_s)$ be the sequence of weights that appear
in the weighted circular graph for $X$. Then
$$
 \itXi(T_{\ast}X)\;
  =\; \{\, (2, a_1),\, \cdots, (2, a_s) \,\}\,.
$$
Consequently, a necessary condition for $T_{\ast}X$ to admit a
splitting type is that the weights that appear in the graph must
be all positive, all zero, or all negative. The augmented
intersection matrix is given by
$$
 Q\; =\; \left[\,
       \begin{array}{cccccc}
        a_1 &   1 &        &        &         &  1 \\
          1 & a_2 &      1 &        &         &    \\
            &   1 &    a_3 &      1 &         &    \\[.6ex]
            &     & \ddots & \ddots & \ddots  &    \\[.6ex]                   
            &     &        &      1 & a_{n-1} &  1 \\
          1 &     &        &        &       1 &  a_n
       \end{array}
              \,\right]\,,
$$
where all the missing entries are $0$. From these data, the
splitting type of $T_{\ast}X$, if exists, can be worked out.

To provide more examples and also for the interest of string theory,
equivariant blowups of $\CP^2$ up to $9$ points that admits a
splitting type are searched out by computer. These include
del Pezzo type and $\frac{1}{2}\,$K3 type surfaces. It turns out
that there are only $8$ of them, $4$ of del Pezzo type and $4$ of
$\frac{1}{2}\,$K3 type. Their splitting types are listed in
{\sc Table 4-1}.

\bigskip

\noindent\hspace{-5em}{\tiny
\begin{tabular}{|c||c|l|c|l|l|} \hline
 \rule{0ex}{3ex} topology of $X$
     & $k$
            & $w=(a_1,\,\cdots,\,a_s)$
                  & $H_2(X;{\tinyBbb Z})$
                  & splitting type of $T_{\ast}X$
                              & Remark    \\[.6ex] \hline\hline 
 \rule{0ex}{4ex}
 ${\tinyBbb C}{\rm P}^2\sharp\,3\,\overline{{\tinyBbb C}{\rm P}^2}$
     & $3$
            & $(-1,-1,-1,-1,-1,-1)$
                  & ${\tinyBbb Z}^{4}$
                  & $\begin{array}{l}
                       (\,(2,4,4,2)\,,\; (-1,-2,-2,-1)\,)
                     \end{array}$
                        & del Pezzo type  \\[.6ex] \hline
 \rule{0ex}{4ex}
 ${\tinyBbb C}{\rm P}^2\sharp\,5\,\overline{{\tinyBbb C}{\rm P}^2}$
     & $5$
            & $(-1,-2,-1,-2,-1,-2,-1,-2)$
                  & ${\tinyBbb Z}^6$
                  & $\begin{array}{l}
                        (\,(2,4,8,6,6,2)\,,\\
                        \hspace{1em}
                        (-2,-3,-6,-4,-4,-1)\,)
                     \end{array}$
                         & del Pezzo type   \\[1.5ex] \hline
 \rule{0ex}{4ex}
 ${\tinyBbb C}{\rm P}^2\sharp\,6\,\overline{{\tinyBbb C}{\rm P}^2}$
     & $6$        
            & $(-1,-2,-2,-1,-2,-2,-1,-2,-2)$
                  & ${\tinyBbb Z}^7$
                  & $\begin{array}{l}
                        (\,(2,4,8,14,8,4,2)\,,\\
                        \hspace{1em}
                        (-2,-3,-6,-11,-6,-3,-2)\,)
                     \end{array}$
                         & del Pezzo type   \\[1.5ex] \hline
 \rule{0ex}{4ex}
 ${\tinyBbb C}{\rm P}^2\sharp\,7\,\overline{{\tinyBbb C}{\rm P}^2}$
     & $7$
            & $(-1,-2,-2,-1,-3,-1,-2,-2,-1,-3)$
                  & ${\tinyBbb Z}^8$
                  & $\begin{array}{l}
                        (\,(2,4,8,14,8,12,6,2)\,,\\
                        \hspace{1em}
                        (-3,-4,-7,-12,-6,-9,-4,-1)\,)
                     \end{array}$
                         & del Pezzo type   \\[1.5ex] \hline
 \rule{0ex}{4ex}
 ${\tinyBbb C}{\rm P}^2\sharp\,9\,\overline{{\tinyBbb C}{\rm P}^2}$
     & $9$
            & $(-1,-2,-2,-2,-1,-4,-1,-2,-2,-2,-1,-4)$
                  & \hspace{1ex}${\tinyBbb Z}^{10}$
                  & $\begin{array}{l}
                      (\,(2,4,8,14,22,10,20,12,6,2)\,, \\
                       \hspace{1em}
                       (-4,-5,-8, -13, -20,-8, -16, -9, -4, -1)\,)
                     \end{array}$
                         & $\frac{1}{2}\,$K3 type  \\[1.5ex] \hline
 \rule{0ex}{4ex}
 ${\tinyBbb C}{\rm P}^2\sharp\,9\,\overline{{\tinyBbb C}{\rm P}^2}$
     & $9$
            & $(-1,-2,-2,-3,-1,-2,-2,-3,-1,-2,-2,-3)$
                  & \hspace{1ex}${\tinyBbb Z}^{10}$
                  & $\begin{array}{l}
                      (\,(2,4,8,14,36,24,14,6,6,2)\,, \\
                       \hspace{1em}
                       (-3,-4,-7,-12,-32,-21,-12,-5,-6,-2)\,)
                     \end{array}$
                         & $\frac{1}{2}\,$K3 type  \\[1.5ex] \hline
 \rule{0ex}{4ex}
 ${\tinyBbb C}{\rm P}^2\sharp\,9\,\overline{{\tinyBbb C}{\rm P}^2}$
     & $9$
            & $(-1,-2,-3,-1,-2,-3,-1,-2,-3,-1,-2,-3)$
                  & \hspace{1ex}${\tinyBbb Z}^{10}$
                  & $\begin{array}{l}
                      (\,(2,4,8,22,16,12,22,12,4,2)\,, \\
                       \hspace{1em}
                       (-3,-4,-7,-20,-14,-10,-19,-10,-3,-2)\,)
                     \end{array}$
                         & $\frac{1}{2}\,$K3 type  \\[1.5ex] \hline
 \rule{0ex}{4ex}
 ${\tinyBbb C}{\rm P}^2\sharp\,9\,\overline{{\tinyBbb C}{\rm P}^2}$
     & $9$
            & $(-1,-3,-1,-3,-1,-3,-1,-3,-1,-3,-1,-3)$
                  & \hspace{1ex}${\tinyBbb Z}^{10}$
                  & $\begin{array}{l}
                      (\,(2,4,12,10,20,12,18,8,8,2)\,, \\
                       \hspace{1em}
                       (-3,-4,-12, -9, -18,-10, -15, -6, -6, -1)\,)
                     \end{array}$
                         & $\frac{1}{2}\,$K3 type  \\[1.5ex] \hline
  \multicolumn{6}{|l|}{\footnotesize
   \rule{0ex}{3ex}
   (1) $k\,$ is the number of points blown up from
                          ${\footnotesizeBbb C}{\rm P}^2$. } \\[.6ex]
  \multicolumn{6}{|l|}{\footnotesize
   \rule{0ex}{0ex}
   (2) $H_2(X;{\footnotesizeBbb Z})$ is generated by
           the first $(s-2)$ divisors in the list $w$.} \\[.6ex]
  \multicolumn{6}{|l|}{\footnotesize
   \rule{0ex}{0ex}
   (3) The line bundles in the splitting type are represented by
       divisors of $X$ as elements in
                $H_2(X;{\footnotesizeBbb Z})$.}  \\[.6ex] \hline
\end{tabular}
} 

\medskip

\begin{center}
 \parbox{11cm}{{\sc Table 4-1}.
 Complete list of $T_{\ast}S$ and its splitting type for toric
 surfaces obtained from $\CP^2$ via equivariant blowups up to $9$
 points.}
\end{center}

\bigskip

The fan for these toric surfaces are indicated in {\sc Figure 4-3}.
 \begin{figure}[htbp]
  \setcaption{{\sc Figure 4-3.}
   \baselineskip 14pt
   The fan for the blowups of $\CP^2$ up to nine points that admit
   a splitting type, together with the weights, are indicated. 
  } 
  \centerline{\psfig{figure=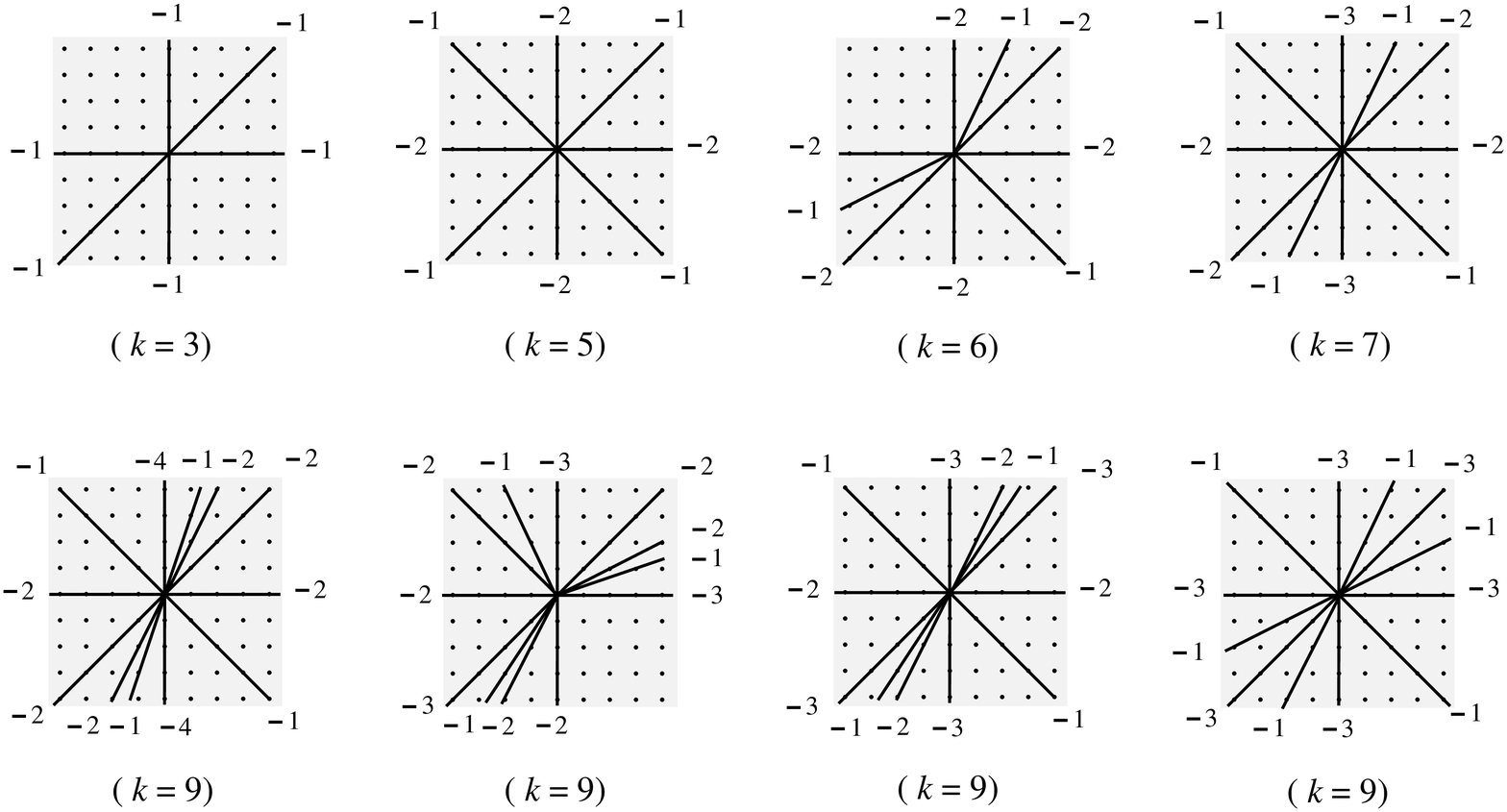,width=13cm,caption=}}
 \end{figure}

This concludes the example.

\noindent\hspace{14cm} $\Box$

\bigskip

\section{Remarks and issues for further study.}

In the previous section, we have illustrated how the data of a
linearized equivariant vector bundle ${\cal E}$ over a toric manifold
$X_{\Sigma}$ is used to determine its splitting type if it exists.
It turns out that these examples are related to the following kind of
exact sequence\footnote{We thank Bong H.\ Lian for drawing our
   attention to this and the reference [Ja]}:
$$
 0\; \longrightarrow\; {\cal O}_{X_{\Sigma}}\;
   \stackrel{\eta}{\longrightarrow}\;
   \oplus_{i=1}^{r+1}\, {\cal O}_{X_{\Sigma}}(D_i)\;
   \longrightarrow {\cal E}\; \longrightarrow\; 0\,,
$$
where $D_i$ are Cartier T-Weil divisors of $X_{\Sigma}$ and $\eta$
is a holomorphic bundle inclusion. For such ${\cal E}$, the system
of splitting numbers $\itXi({\cal E})$ may be obtained directly from
this exact sequence.

\bigskip

\noindent
{\bf Example 5.1
[equivariant vector bundles of rank $n$ over ${\Bbb C}{\rm P}^n$].}
Let $[z_0:\,\cdots,\,z_n]$ be the homogeneous coordinates of the
projective space ${\Bbb C}{\rm P}^n$. Recall from [Ka2] that an
indecomposable equivariant vector bundle ${\cal E}$ of rank $n$
over ${\Bbb C}{\rm P}^n$ is isomorphic to
either $E\otimes{\cal O}_{{\scriptsizeBbb C}{\rm P}^n}(d)$ or
$E^{\ast}\otimes{\cal O}_{{\scriptsizeBbb C}{\rm P}^n}(d)$ for
some integer $d$ where $E$ is the equivariant vector bundle defined
by the exact sequence
$$
 0\; \longrightarrow\;
 {\cal O}_{{\scriptsizeBbb C}{\rm P}^n}\;
 \stackrel{\eta}{\longrightarrow}\;
 \oplus_{i=0}^n\,{\cal O}_{{\scriptsizeBbb C}{\rm P}^n}(m_i)\;
 \longrightarrow\; E\; \longrightarrow\; 0\,,
$$
where $m_i$ are positive integers and $\eta$ sends $1$ to
$(z_0^{m_0}, \,\cdots,\,z_n^{m_n})$. Since the $(n+1)n/2$ many
invariant ${\Bbb C}{\rm P}^1$ in ${\Bbb C}{\rm P}^n$ are given by
$$
 V_{ij}\;=\; \{[\,0,\,\cdots,\,0,\,z_i,\,0,\,\cdots,\,0,\,z_j,\,
     0,\,\cdots,\,0\,]\,|\, (z_i,z_j)\in {\Bbb C}^2-\{(0,0)\}\}\,,
$$
$0\le i<j\le n$, the above exact sequence, when restricted to
$V_{ij}$, reduces to
$$
 \begin{array}{cccccccccl}
   0  & \longrightarrow   & {\cal O}_{V_{ij}}   & \longrightarrow
      & \oplus_{k=1}^{r+1}{\cal O}_{V_{ij}}(m_k)
      & \longrightarrow   & E|_{V_{ij}}
      & \longrightarrow   & 0                    \\[1ex]
  & & 1  & \longmapsto
      & (\,0,\,\cdots,\,0,\, z_i^{m_i},\, 0,\, \cdots,\,0,\,
                                 z_j^{m_j},\, 0,\,\cdots,\,0\,)
      &  &  & & &.
 \end{array}
$$
It follows from the multiplicativity of total Chern class that
$$
 E|_{V_{ij}}\;
 \simeq\; {\cal O}_{{\scriptsizeBbb C}{\rm P}^1}(m_i+m_j)\,
           \oplus\,
           {\cal O}_{{\scriptsizeBbb C}{\rm P}^1}(m_0)\,
            \oplus\, \cdots\, \oplus\,
           \widehat{ {\cal O}_{{\scriptsizeBbb C}{\rm P}^1}(m_i) }\,
           \oplus\, \cdots\, \oplus\,
           \widehat{ {\cal O}_{{\scriptsizeBbb C}{\rm P}^1}(m_j) }\,
           \oplus\, \cdots\, \oplus\,
           {\cal O}_{{\scriptsizeBbb C}{\rm P}^1}(m_n)
$$
and, hence, the system of splitting numbers of $E$ is
$$
 \itXi(E)\;
 =\; \{\, (m_i+m_j, m_0,\,\cdots,\, \widehat{m_i},\,\cdots,\,
        \widehat{m_j},\,\cdots,\,m_n)\,|\, 0\le i<j\le n \,\}\,,
$$
where terms with $\widehat{\;\;}$ are deleted.
The augmented matrix in this case is
$$
 Q\;=\; \left[ \begin{array}{ccc}
                  1       & \cdots  & 1 \\
                  \vdots  & \vdots  & \vdots \\
                  1       & \cdots  & 1 \\
               \end{array}      \right]\;.
$$
Without loss of generality, one may assume that
$0<m_0\le\cdots\le m_n$, then $m_i< m_n+m_{n-1}$ for all $i$.
Consequently, with the notation from Sec.\ 3, for $QX=R^{\prime}$
to have a solution, one must have $m_i+m_j=m_n+m_{n-1}$ for all
$i<j$, which implies that $m_0=\cdots=m_n$. One concludes therefore

\bigskip

\noindent
{\bf Corollary.} {\it For the indecomposable equivaraint rank $n$
 bundle ${\cal E}=E\otimes{\cal O}_{{\scriptsizeBbb C}{\rm P}^n}(d)$
 over $\CP^n$ to admit a splitting type, one must have
 $m_0=\cdots=m_n$. In this case, the splitting type of ${\cal E}$
 is unique and is given by
 $(\,{\cal O}_{{\scriptsizeBbb C}{\rm P}^n}(2m_0+d),\,\cdots,\,
     {\cal O}_{{\scriptsizeBbb C}{\rm P}^n}(m_0+d)\,)$.
} 

\bigskip

This generalizes Example 4.1. Note also that the case
$m_0=\cdots=m_n=1$ with $d=0$ corresponds to
$T_{\ast}{\Bbb C}{\rm P}^n$ and the above discussion double-checks
part of Example 4.2.

\noindent\hspace{14cm} $\Box$

\bigskip

One can generalize this example slightly to toric manifolds as
follows. First, let us state a lemma, whose proof is straightforward.

\bigskip

\noindent
{\bf Lemma 5.2} {\it
 Given an exact sequence 
 $$
  0\; \longrightarrow\;
    {\cal O}_{{\scriptsizeBbb C}{\rm P}^1}\;
    \stackrel{\eta}{\longrightarrow}\;
    {\cal O}_{{\scriptsizeBbb C}{\rm P}^1}\;
      \oplus\left(
       \oplus_{i=1}^{r}
         {\cal O}_{{\scriptsizeBbb C}{\rm P}^1}(m_i)\right)\;
    \longrightarrow\; E\; \longrightarrow\; 0\,,
 $$
 where $\eta(1)=(s_0, s_1,\cdots, s_r)$, such that $s_0$ is
 non-zero. Then
 $E\simeq \oplus_{i=1}^r{\cal O}_{{\scriptsizeBbb C}{\rm P}^1}(m_i)$.
} 

\bigskip

Given a toric $n$-fold $X_{\Sigma}$, consider now the exact sequence
$$
 0\; \longrightarrow\; {\cal O}_{X_{\Sigma}}\;
   \stackrel{\eta}{\longrightarrow}\;
   \oplus_{i=1}^{r+1}\, {\cal O}_{X_{\Sigma}}(D_i)\;
   \longrightarrow {\cal E}\; \longrightarrow\; 0\,.
$$
The restriction of the sequence to $V(\tau)$, $\tau\in\Sigma(n-1)$,
is given by
$$
 0\; \longrightarrow\;
   {\cal O}_{{\scriptsizeBbb C}{\rm P}^1}\;
   \stackrel{\eta}{\longrightarrow}\;
      \oplus_{i=1}^{r+1}
        {\cal O}_{{\scriptsizeBbb C}{\rm P}^1}(D_i\cdot V(\tau))\;
   \longrightarrow\; E\; \longrightarrow\; 0\,.
$$
If this exact sequence is of the kind in Example 5.1 or Lemma 5.2
for all $\tau\in\Sigma(n-1)$, then $\itXi({\cal E})$ can be readily
obtained and the splitting type of such ${\cal E}$, if exists, can
then be determined. Inspired from Example 5.1, to realize this, recall
the Cox homogeneous coordinates of $X_{\Sigma}$ from [Co]
(cf.\ Sec.\ 1):
let $a=|\Sigma(1)|$, then $X_{\Sigma}$ can be realized as a quotient
$X_{\Sigma}=({\Bbb C}^{\Sigma(1)}-Z(\Sigma))/\mbox{\raisebox{-.4ex}{$G$}}$.
Let $(z_1,\,\cdots,\,z_a)$ be the standard coordinates of ${\Bbb C}^a$
and $\tau=[v_{j_1},\,\cdots,\,v_{j_{n-1}}]\in\Sigma(n-1)$, then
$V(\tau)$ can be realized as the quotient of the coordinate subspace:
$V(\tau)=\{z_{j_1}=\,\cdots\,=z_{j_{n-1}}\}/\mbox{\raisebox{-.4ex}{$G$}}$.
Furthermore, if $\tau=\sigma_1\cap\sigma_2$, where
$$
 \sigma_1=[v_{j_1},\,\cdots,\,v_{j_{n-1}}, v_{j_n}]
 \hspace{1em}\mbox{and}\hspace{1em}
 \sigma_2=[v_{j_1},\,\cdots,\,v_{j_{n-1}}, v_{j^{\prime}_n}]\,,
$$
then $[z_{j_n}:z_{j^{\prime}_n}]$ serves as a homogeneous coordinates
for $V(\tau)\simeq\CP^1$. For all other $i$, 
$\{z_i=0\}\cap\{z_{j_1}=\,\cdots\,=z_{j_{n-1}}\}$ lies in the
exceptional subset $Z(\Sigma)$ and, hence, $z_i$ as an element in
the homogeneous coordinate ring ${\Bbb C}[z_1,\,\cdots,\,z_a]$,
graded by the Chow group $A_{n-1}(X_{\Sigma})$, descends to a
non-zero section in
${\cal O}_{X_{\Sigma}}(D(v_i))|_{V(\tau)}
                    \simeq{\cal O}_{{\scriptsizeBbb C}{\rm P}^1}$.
In general, since
${\cal O}_{X_{\Sigma}}(D_1)\otimes{\cal O}_{X_{\Sigma}}(D_2)
  ={\cal O}_{X_{\Sigma}}(D_1+D_2)$
for any Cartier T-Weil divisor $D_1$, $D_2$, for any monomial
$\prod_{k}z_{j_k}^{\alpha_k}$ with
$j_k\notin\{j_1,\,\cdots,\,j_{n-1}, j_n, j^{\prime}_n\}$ and
$\alpha_k$ positive integers, $\prod_{k}z_{j_k}^{\alpha_k}$ descends
to a non-zero section in
${\cal O}_{X_{\Sigma}}(\sum_k\alpha_kD(v_{j_k}))|_{V(\tau)}
                    \simeq{\cal O}_{{\scriptsizeBbb C}{\rm P}^1}$.
This fact provides us with a guideline for defining $\eta$ so that
exact sequences as in Lemma 5.2 can appear when restricted to
invariant $\CP^1$'s in $X_{\Sigma}$.
Such examples can be constructed plenty. Let us give an example below
to illustrate the idea.

\bigskip

\noindent
{\bf Example 5.3 [simple rank 3 bundle over Hirzebruch surface].}
Let $X_{\Sigma}={\Bbb F}_a$ be a Hirzebruch surface
(cf. Example 4.2 (b)).
Consider the rank $3$ bundle ${\cal E}(m_1, m_2, m_3, m_4)$
over ${\Bbb F}_a$ defined by the exact sequence
$$
 \begin{array}{cccccccccl}
   0  & \longrightarrow   & {\cal O}_{{\scriptsizeBbb F}_a}
      & \longrightarrow
      & \oplus_{k=1}^4{\cal O}_{{\scriptsizeBbb F}_a}(m_k D_{v_k})
      & \longrightarrow   & {\cal E}(m_1, m_2, m_3, m_4)
      & \longrightarrow   & 0                    \\[1ex]
  & & 1  & \longmapsto
      & (\,z_1^{m_1},\, z_2^{m_2},\, z_3^{m_3},\, z_4^{m_4})
      &  &  & & &,
 \end{array}
$$
where $m_i$ are positive integers. From Lemma 5.2 and the discussions
above, one concludes that the system of splitting numbers of
${\cal E}(m_1, m_2, m_3, m_4)$ is given by
$$
 \itXi({\cal E})\;
 =\;  \{\, (m_1, m_2, m_4), (m_1, m_2, m_3), (m_2, m_3, m_4),
              (m_1, m_3, m_4)\,\}\,.
$$
Recall the augmented intersection matrix $Q$ for ${\Bbb F}_a$ from
Example 4.2 (b). Through a tedious but straightforward algebra, one
can show that the only case when ${\cal E}(m_1, m_2, m_3, m_4)$
admits a splitting type is when $a=0$ (i.e.\ $X=\CP^1\times\CP^1$)
with $m_1=m_3$ and $m_2=m_4$. In this case, the splitting type is
unique and is given by
$$
 {\cal O}_{{\scriptsizeBbb C}{\rm P}^1
     \times{\scriptsizeBbb C}{\rm P}^1}(m_2, m_1)\,
  \oplus\, {\cal O}_{{\scriptsizeBbb C}{\rm P}^1
              \times{\scriptsizeBbb C}{\rm P}^1}(m_1, m_2)\,
  \oplus\, {\cal O}_{{\scriptsizeBbb C}{\rm P}^1
                      \times{\scriptsizeBbb C}{\rm P}^1}(m_1, m_2)\,,
$$
where we idetify $\Pic(X)$ with
$H_2(X,{\Bbb Z})\simeq{\Bbb Z}\oplus{\Bbb Z}$, with generators
$$
 {\cal O}_{{\scriptsizeBbb C}{\rm P}^1
   \times{\scriptsizeBbb C}{\rm P}^1}(1,0)\;
                \longmapsto\;[\CP^1\times\ast]\,, \hspace{2em}
{\cal O}_{{\scriptsizeBbb C}{\rm P}^1
   \times{\scriptsizeBbb C}{\rm P}^1}(0,1)\;
                \longmapsto\;[\ast\times\CP^1]\,.
$$

\noindent\hspace{14cm} $\Box$

\bigskip

For more general $\eta$, the restriction of the exact sequence over
$X_{\Sigma}$ to each invariant ${\Bbb C}{\rm P}^1$ in $X_{\Sigma}$
leads to an exact sequence of the form
$$
 0\; \longrightarrow\;
   {\cal O}_{{\scriptsizeBbb C}{\rm P}^1}\;
   \stackrel{\eta}{\longrightarrow}\;
   \oplus_{i=1}^{r+1}\, {\cal O}_{{\scriptsizeBbb C}{\rm P}^1}(m_i)\;
   \longrightarrow {\cal E}_{{\scriptsizeBbb C}{\rm P}^1}\;
   \longrightarrow\; 0\,.
$$
A complete study of how $\eta$ determines the splitting of
${\cal E}_{{\scriptsizeBbb C}{\rm P}^1}$ as a direct sum of line
bundles requires more work\footnote{
 We thank Jason Starr for discussions on this and the references
 [Bri], [E-VV1], and [E-VV2].}.

\bigskip

We conclude the discussion of splitting types here and leave its
further study and applications for another work.

\newpage

\begin{flushleft}
{\large\bf Appendix. The computer code.}
\end{flushleft}
The computer code in Mathematica that carries out the computation
in Example 4.2 (c) is attached below for reference,

\bigskip

{\scriptsize
\baselineskip 9.6pt
\begin{verbatim}

(* This is a code in Mathematica. *)
(* The purpose of this code is to sort out and compute the splitting type of the tangent bundle
   of toric surfaces. The result of computation is written to the file 'ma-result.txt'. *)
(* Subroutines enclosed: BlowUp, BlowUpN, GenerateMatrix, (MAIN) SplittingType *)

(* Definition of the function 'BlowUp'. *)
(* 'BlowUp[weightlist]' generates the list of weights on the circular weighted graph obtained by
    equivariant blowup at a $T_N$-fixed point of a toric surface represented by 'weightlist'.
   Date of completion: 10/15/1999. Test: Tested correct. Date of last revision: 10/16/1999.
 *)

  BlowUp[ weightlist_ ] :=
    Module[ { a1, a2, b, b1, b2, list, list1, list2, list3, m1, m2, newlist },

            m1=Length[weightlist];

            list1[i_] := ReplacePart[ weightlist,
                           { weightlist[[i]]-1, -1, weightlist[[i+1]]-1 }, i+1 ];
            list2[i_] := Delete[ list1[i], i];
            list3[i_] := Flatten[ list2[i] ] ;

            list=ReplacePart[ weightlist, {-1, weightlist[[1]]-1 }, 1 ];
            list=ReplacePart[ list, list[[m1]]-1, m1];
            list={ Flatten[list] };

            newlist=Join[ list, Table[ list3[i] , {i, 1, m1-1}] ];
            newlist=Union[newlist];
            m2=Length[newlist];

            Do[
                a1=newlist[[i]];
                a2=Reverse[a1];
                b1=Table[ RotateRight[a1, i], {i, 1, m1+1} ];
                b2=Table[ RotateRight[a2, i], {i, 1, m1+1} ];
                b=Union[b1, b2];
                newlist=Union[{a1}, Complement[newlist, b] ];
                If[ m2>Length[newlist],
                    Return[newlist]
                  ],
               {i, 1, m2}
              ];

            Return[newlist]
          ]


(* Definition of the function 'BlowUpN'. *)
(* 'BlowUpN[weightlist, n]' generates the list of weights on the circular weighted graph obtained
    by consecutive equivariant blowup at a $T_N$-fixed points of a toric surface 'n' times, starting
    from the one represented by 'weightlist'.
   Date of completion: 10/15/1999. Test: Tested correct. Date of last revision: 10/16/1999.
 *)

  BlowUpN[ weight_, n_ ] :=
    Module[ { m, newlist, oldlist, totallist},

            totallist={weight};
            oldlist={weight};
            newlist={};

            Do[
                m=Length[oldlist];
                Do[
                    newlist=Union[ newlist, BlowUp[ oldlist[[j]] ] ],
                   {j, 1, m}
                  ];
                oldlist=newlist;
                totallist=Join[ totallist, newlist],
               {i, 1, n}
              ];
            totallist=Union[totallist];
            Return[totallist];
          ]


(* Definition of the function 'GenerateMatrix'. *)
(* 'GenerateMatrix[weightlist]' generates a matrix following the rule discussed in the paper on
    splitting types of equivariant vector bundle on toric manifolds .
   Date of completion: 10/15/1999. Test: Tested correct. Date of last revision: 10/15/1999.
 *)

  GenerateMatrix[ weightlist_ ] :=
    Module[ { listfirst, listlast, list1, list2, m, newlist, v  },

            m=Length[weightlist];
            v=Table[ 0, {i, 1, m-2} ];

            list1[i_]:=ReplacePart[ v, { 1, weightlist[[i]], 1 }, i-1 ];
            list2[i_]:=Flatten[ list1[i] ];

            listfirst=Flatten[ ReplacePart[ v, {1, weightlist[[1]], 1 }, 1 ] ];
            listfirst={ RotateLeft[listfirst, 1] };

            listlast=Flatten[ ReplacePart[ v, {1, weightlist[[m]], 1 }, 1 ] ];
            listlast={ RotateLeft[listlast, 2] };

            newlist=Join[ listfirst, Table[ list2[i], {i, 2, m-1} ], listlast ];

            Return[newlist]
          ]



(* MAIN ROUTINE *)

(* Definition of the function 'SplittingType'. *)
(* 'SplittingType[weight, n]' sorts out from all the toric surfaces that arise from equivariant
    blowups up to 'n' times of the toric surface whose associated weighted circular graph is given
    by 'weight' those that admit a splitting type and computes their splitting types.
   Date of completion: 10/15/1999. Test: Tested correct. Date of last revision: 10/17/1999.
 *)

 SplittingType[weight_, n_] :=
   Module[ { b, m, matrix, t, totallist, t1, x1, x2 },

           totallist=BlowUpN[weight, n];
           m=Length[totallist];

           Do[
               t=totallist[[i]];
               t1=Union[t];
               If[ Complement[t1,{0}]===t1,
                   matrix=GenerateMatrix[t];
                   b=Table[2, {j, 1, Length[t]}];
                   x1=LinearSolve[matrix, b1];
                   If[ Length[x1]>=3,
                       x2=LinearSolve[matrix, t];
                       PutAppend[i, "ma-result.txt"];
                       PutAppend[t, "ma-result.txt"];
                       PutAppend[x1, "ma-result.txt"];
                       PutAppend[x2, "ma-result.txt"]
                     ]
                 ],
              {i, 1, m}
             ];
         ]



(* Case of study *)

  DeleteFile["ma-result.txt"];
  SplittingType[{1, 1, 1}, 9];
  

\end{verbatim}

} 

\bigskip

\newpage


\end{document}